\documentclass{article}
\usepackage[margin=3cm]{geometry}
\usepackage{amsmath}
\usepackage{graphicx}
\usepackage{enumerate}
\usepackage[numbers]{natbib}
\usepackage{url}

\usepackage{amssymb}
\usepackage{xcolor}
\usepackage{booktabs}
\usepackage{multirow}
\usepackage{algorithm}
\usepackage{algpseudocode}
\def\bmu{\boldsymbol{\mu}}

\def\btheta{\boldsymbol{\theta}}

\def\bSigma{\boldsymbol{\Sigma}}
\def\bdelta{\boldsymbol{\delta}}
\def\bDelta{\boldsymbol{\Delta}}

\def\b0{\mathbf{0}}
\def\bx{\mathbf{x}}
\def\bA{\mathbf{A}}
\def\bB{\mathbf{B}}
\def\bD{\mathbf{D}}
\def\be{\mathbf{e}}
\def\bg{\mathbf{g}}
\def\bp{\mathbf{p}}
\def\bq{\mathbf{q}}
\def\bO{\mathbf{O}}
\def\bP{\mathbf{P}}
\def\bI{\mathbf{I}}
\def\bG{\mathbf{G}}

\def\bC{\mathbf{C}}
\def\bR{\mathbf{R}}

\def\ba{\mathbf{a}}
\def\bj{\mathbf{j}}
\def\br{\mathbf{r}}
\def\bu{\mathbf{u}}
\def\bv{\mathbf{v}}

\def\bx{\mathbf{x}}
\def\by{\mathbf{y}}

\def\bX{\mathbf{X}}

\def\cI{{\cal I}}
\def\cG{{\cal G}}
\def\cJ{{\cal J}}
\def\cN{{\cal N}}
\def\cP{{\cal P}}
\def\DsKL{D_{\mathrm{sKL}}}

\def\bumu{\mathbf{u}_{\boldsymbol{\mu}}}
\def\buSigma{\mathbf{u}_{\boldsymbol{\Sigma}}}
\def\bvmu{\mathbf{v}_{\boldsymbol{\mu}}}
\def\bvSigma{\mathbf{v}_{\boldsymbol{\Sigma}}}

\newcommand{\deb}[1]{\frac{d}{d #1}}
\newcommand{\de}[2]{\frac{d #1}{d #2}}
\newcommand{\den}[3]{\frac{d^{#3} #1}{d #2^{#3}}}

\newcommand{\pdeb}[1]{\frac{\partial}{\partial #1}}
\newcommand{\trop}[1]{\mathrm{tr}\left(#1\right)}

\DeclareMathOperator{\acosh}{acosh}
\DeclareMathOperator{\sign}{sign}
\DeclareMathOperator{\tr}{tr}

\newcommand{\tabref}[1]{Table~\ref{#1}}
\newcommand{\figref}[1]{Fig.~\ref{#1}}
\newcommand{\footremember}[2]{%
    \footnote{#2}
    \newcounter{#1}
    \setcounter{#1}{\value{footnote}}%
}
\newcommand{\footrecall}[1]{%
    \footnotemark[\value{#1}]%
}

\title{\bf The Fisher Geometry and Geodesics of the Multivariate Normals, without Differential Geometry}
\author{%
    Brodie A. J. Lawson\footremember{CDS}{Centre for Data Science, Queensland University of Technology, Australia}\footremember{QUT}{School of Mathematical Sciences, Queensland University of Technology, Australia}\footremember{ARCp}{ARC Centre of Excellence for Plant Success in Nature and Agriculture, Australia}%
    \and Kevin Burrage\footrecall{QUT} \footrecall{ARCp} \footremember{Oxford}{Department of Computer Science, University of Oxford, United Kingdom (Visiting Professorship)}%
    \and Kerrie Mengersen\footrecall{CDS} \footrecall{QUT}%
    \and Rodrigo Weber dos Santos\footremember{UFJF}{Graduate Program in Computational Modeling, Universidade Federal de Juiz de Fora, Brazil}%
    }
    \date{}
\begin{document}
  
\maketitle

\begin{abstract}
Choosing the Fisher information as the metric tensor for a Riemannian manifold provides a powerful yet fundamental way to understand statistical distribution families. Distances along this manifold become a compelling measure of statistical distance, and paths of shorter distance improve sampling techniques that leverage a sequence of distributions in their operation. Unfortunately, even for a distribution as generally tractable as the multivariate normal distribution, this information geometry proves unwieldy enough that closed-form solutions for shortest-distance paths or their lengths remain unavailable outside of limited special cases. In this review we present for general statisticians the most practical aspects of the Fisher geometry for this fundamental distribution family. Rather than a differential geometric treatment, we use an intuitive understanding of the covariance-induced curvature of this manifold to unify the special cases with known closed-form solution and review approximate solutions for the general case. We also use the multivariate normal information geometry to better understand the paths or distances commonly used in statistics (annealing, Wasserstein). Given the unavailability of a general solution, we also discuss the methods used for numerically obtaining geodesics in the space of multivariate normals, identifying remaining challenges and suggesting methodological improvements.
\end{abstract}

\noindent%
\vfill

\newpage

\section{Introduction}

The concept of {\it statistical distance}, some measure of discrepancy between two distributions, is ubiquitous in statistics. Classical tests such as the Cram{\'{e}}r--von Mises criterion or Anderson--Darling test use a distance between (cumulative) density functions to identify whether given data is expected to have come from a given distribution~\citep{DasGupta2008}. The expectation maximisation algorithm, and variational approaches to Bayesian inference, operate by minimising an information theoretic discrepancy measure --- the Kullback--Leibler (KL) divergence~\citep{Bishop2006}. The Wasserstein distance (also known as the Kantrovich--Rubenstein metric or earth-mover's distance) describes in a sense the amount of work required to physically reshape one distribution into another, a concept that brings with it both appealing mathematics~\citep{Li2023}, and demonstrated potential in statistics and machine learning~\citep{Arjovsky2017, Srivastava2018, Panaretos2019}. These are only a few examples of the many statistical discrepancy measures that see use~(for example see also \citep{Jeffreys1946, vanErven2014}).

Spoiled for choice in such a manner, it is often desirable to select a discrepancy that is mathematically specified by a distance metric, such that we may understand our measure of statistical distance through analogs to the familiar Euclidean distance (satisfying key properties such as symmetry, positivity, and the triangle inequality). If we also choose this metric such that our notion of distance is tied to statistical information, then lengths of paths through distribution space become measures of the information required to traverse that path, and shortest-distance paths become ``paths of least statistical resistance''. This is of practical concern, as can be seen within sampling approaches using distribution sequences such as sequential Monte Carlo~\citep{Sim2012} or parallel tempering~\citep{Syed2021}, as well as normalising constant (ratio) estimators such as annealed importance sampling~\citep{Gelman1998, Grosse2013}.

The Fisher information matrix serves as just such a (Riemannian) metric. Use of this metric began with Rao in 1945, who wile presenting the well-known Cram{\'{e}}r--Rao bound in statistical estimation also discussed the idea of members of a distribution family sitting on what is now known as a ``statistical manifold''~\citep{Rao1945}. The geometry of this Riemannian manifold is warped by the information already present in the distribution, in a fashion similar to the way Riemannian manifolds in astrophysics describe the gravitational warping of space-time~\cite[Ch. 10]{Landau1971}. Especially after their modernisation by \citet{Amari1985}, these ideas of {\it information geometry} have often concentrated on the theoretical properties of the manifolds associated with statistical models or distribution families~\citep{Nielsen2020}, but have also seen use in generating new statistical algorithms~\citep{Banerjee2005}, explaining convergence behaviour~\citep{Ikeda2006} and making learning algorithms parameterisation-invariant~\citep{Martens2020}. As a means of quantifying statistical distance, information geometry has a natural application to hypothesis testing~\citep{Calvo1991}, and is also intrinsically linked to the question of identifiability of parameters in models of complex systems~\citep{Transtrum2011, Sharp2022}.

For a given family of distributions, information geometry provides an information-aware measure of distance termed the Fisher--Rao distance (or sometimes the Rao distance) that behaves like a distance, unlike say the KL divergence~\citep{vanErven2014}. However, owing to the way Riemannian geometry defines distance in an infinitesimal fashion, finding the Fisher--Rao distance between two points on the statistical manifold requires first finding the shortest-distance path along the manifold between those two points. Known as geodesics, such paths generalise straight lines to a non-Euclidean geometry, but are challenging to generate. This, together with the general requirement of familiarity with differential geometry, means that information geometry remains somewhat unwieldy.

In this work, we provide an introduction to, and practical understanding of, information geometry for general statisticians. We do so by avoiding as many differential geometric concepts and as much of the corresponding notation as possible, and by concentrating on a fundamental distribution family, the multivariate normals (MVNs). What may at first appear to be a simplistic choice reveals a rich topic area of ongoing research that we review, with applications for example to segmentation of diffusion tensor imaging data~\citep{Han2014} and wireless signal processing~\citep{Pilte2016,Du2020}. We provide an intuitive way of understanding this statistical manifold's geometry, and the form taken by its geodesic paths. This approach allows us to unify the existing cases for which a closed-form expression for the geodesic is known, and informs approximations for use when a closed-form solution is not available. We also discuss how other commonly-used paths through distribution space (annealing/geometric, Wasserstein-optimal transport) compare to these information-optimal paths in the space of MVNs, providing a deeper understanding of their function and expected performance.

Given the need to numerically find the geodesic in the majority of cases where a closed-form solution is unavailable, we also thoroughly discuss the shooting methods used for this purpose. As part of our review, we link recent developments in approximately-geodesic curves to shooting method performance, by demonstrating the potential of the former as an initialisation. We conclude with a discussion of how the concepts we discuss extend beyond the MVN family, and the unsolved problem of a closed form for the geodesic and Fisher--Rao distance between arbitrary MVNs.

\section{The Fisher Geometry of the Normal Distribution}

\subsection{The Riemannian Metric}

Consider a family of distributions parameterised by $\btheta$. For the MVN family we consider here, the typical choice of parameters is the mean vector, $\bmu$, and covariance matrix, $\bSigma$. As such, the parameter space for $d$-variate MVNs could be taken as $\btheta \in \mathbb{R}^{d + d(d+1)/2}$, with the vector $\btheta$ housing the $d$ components of $\bmu$ and the $d(d+1)/2$ unique elements of the symmetric $\bSigma$. Owing to the requirement that $\bSigma$ be a valid (positive definite) covariance matrix, the parameter space for the family of $d$-variate normals in fact occupies some subspace of $\mathbb{R}^{d + d(d+1)/2}$.

Asked for the distance between two points in the parameter space for $\btheta$, we could say calculate the Euclidean norm of the difference in their associated vectors, $\|\btheta_2 - \btheta_1\|_2$. This would completely ignore the statistical information represented within $\btheta$, through its constituent components $\bmu$ and $\bSigma$ and what they imply about the MVN-distributed random variables. Information geometry rectifies this, by treating the parameter space instead as a curved manifold that explicitly acknowledges this statistical context.

As stated, our goal here is to avoid a differential geometric treatment of the statistical manifold associated with MVNs, which has already been comprehensively provided~\citep{Calvo1991, Skovgaard1984, Nielsen2023}. However, in describing the (Riemannian) geometry of this manifold, the metric tensor is an inescapable concept. Commonly denoted $g_{ij}$, the metric tensor defines the warping of space by defining the distance $ds$ that is traversed when the co-ordinates are infinitesimally incremented by some adjustment $d \btheta$,
\begin{equation}
\label{metric}
ds^2 = \sum_i \sum_j g_{ij}\, d\theta_i d\theta_j = d\btheta^T \bG d\btheta,
\end{equation}
with $\theta_i$ denoting the $i$-th co-ordinate of $\btheta$, and $\bG$ the matrix with elements $g_{ij}$. When $\bG$ is the identity matrix, we have $ds^2 = d\btheta^T d\btheta = \|d \btheta\|_2^2$ and so recover the familiar Euclidean distance for a non-warped space. On the other hand, where $\bG$ is itself a function of the parameters, as it is for the MVN statistical manifold, the distance associated with $d \btheta$ depends also on the current position on this manifold.

Information geometry typically imposes this curvature on the parameters of a family of distributions by choosing $\bG$ to be the Fisher information matrix associated with that distribution family~\citep{Amari1985}. For a family with density functions $f(\bx;\btheta)$, the Fisher information matrix $\cI$ has elements defined by
\[
\cI_{ij} = E_{\bX}\left[ \pdeb{\theta_i} \log f(\bx; \btheta) \, \pdeb{\theta_j} \log f(\bx; \btheta) \right].
\]
For $d$-variate MVNs with log-density
\[
\log f(\bx;\bmu, \bSigma) = -\frac{d}{2} \log 2\pi -\frac{1}{2} \log \det \bSigma - \frac{1}{2} (\bx - \bmu)^T \bSigma^{-1} (\bx - \bmu),
\]
we may show that the distance metric takes the form~\citep{Herntier2022}
\begin{equation}
    \label{MVNmetric}
ds^2 = d\bmu^T \bSigma^{-1} d\bmu + \frac{1}{2}\tr\Bigl[ (\bSigma^{-1} d\bSigma )^2 \Bigr] = \left\|\bSigma^{-1/2}d\bmu\right\|^2 + \frac{1}{2} \left\|\bSigma^{-1/2} d\bSigma \bSigma^{-1/2}\right\|_F^2,
\end{equation}
where $\|\cdot\|_F$ denotes the Frobenius norm. Equation~\eqref{MVNmetric} completely defines the geometry of the Fisher manifold for MVNs, but only for two points infinitesimally close together (separated by $\bmu$ and $\bSigma$), and it is not immediately obvious what minimum distance paths between finitely-separated points should look like. Before we discuss such paths, however, it is beneficial to gain a high-level understanding of how the warping of space defined by equation~\eqref{MVNmetric} acts.

The first thing we may observe is that for any point on the manifold $\btheta = (\bmu, \bSigma)$, the distance is affected only by $\bSigma$ (which appears only as the precision, $\bSigma^{-1}$). This makes good statistical sense, in that different values of $\bmu$ do not inherently correspond to differing amounts of information, only to information pointing to different locations in the sample space for the associated random variable. The precision matrix, however, describes how well $\bmu$ is specified, and thus has direct effect on how the distribution will respond to additional information. That is, the greater the precision, the more information that is required to obtain a shift in the distribution's parameters. This information cost is composed of wholly separate contributions due to the movement in $\bmu$-spaec and the movement in $\bSigma$-space, owing to the block-diagonal structure of the Fisher information matrix in this case~\citep{Malago2015}.

The second piece of intuition comes from the norm-based definition for $ds^2$ in equation~\eqref{MVNmetric}. Considering the property for an MVN-distributed random variable $\bx \sim \cN(\bmu, \bSigma)$ that the affine-transformed variable $\by = \bP \bx + \br$ is distributed $\by \sim \cN(\bP \bmu + \br, \bP \bSigma \bP^T)$, we describe affine transformations of an MVN's parameters as those of form
\begin{equation}
    \label{affine_transform}
\bmu \mapsto \bP \bmu + \br, \qquad \qquad \bSigma \mapsto \bP \bSigma \bP^T, \qquad \qquad \det(\bP) > 0.
\end{equation}
We then see through the second expression in equation~\eqref{MVNmetric} that the elements of $d\btheta = (d\bmu, d\bSigma)$ do in fact add as a sum of squares reminiscent of Euclidean distance, but only after the effect of the current precision is removed by applying the affine transformation $\bP = \bSigma^{-1/2}$, $\br = \b0$ that achieves an identity covariance (known as a sphering or whitening transformation). Indeed, the well-known Mahalanobis distance~\citep{Mahalanobis1936} may be thought of as the Euclidean distance after application of a whitening transformation, and the first term in~\eqref{MVNmetric} is precisely the (squared) Mahalanobis distance between two means separated by $d\bmu$. This crystallises the property of the MVN statistical manifold that shifts in the mean in strongly-specified directions cost more information than shifts in weakly-specified directions.

It can also be shown via direct substitution that the metric tensor for MVNs~\eqref{MVNmetric} is invariant to further affine transformations of the form~\eqref{affine_transform}. Interpreted through the lens of the norm-based definition, this invariance may be intuitively understood as pure translations of the mean (through $\br$) not affecting $d\bmu$, and the sphering transformation applied in~\eqref{affine_transform} removing any scaling effects imparted by $\bP$ (rotation may still occur). This invariance of the metric is an important property, as for a given pair of points, we may first transform these points so that one lies at the origin, $(\bmu, \bSigma) = (\b0, \bI)$, without affecting the lengths or intrinsic shapes of these paths~\citep{Imai2011} between them.

\subsection{Geodesics in Multivariate Normal Distribution Space}

The distance metric~\eqref{MVNmetric} directly specifies infinitesimal distances on the statistical manifold associated with MVNs. To define {\it finite} distances on this manifold, these infinitesimal increments must be integrated along the path of minimal distance that connects them --- the geodesic. In Euclidean space, the geodesics are straight lines, and thus the Euclidean metric functions for both infinitesimal and finite distances. For MVNs, however, the geodesics are nontrivial and unavailable in closed form. When found, however, their length gives the Fisher--Rao distance between a pair of MVNs, a statistically-informed measure of their separation. For very close members of a distribution family (MVN or otherwise), half of the square of the Fisher--Rao distance approximates the KL divergence~\citep{Calin2014}. For more distant distributions, the two become distinct.

\begin{figure}
    \centering
    \includegraphics[width=0.48\textwidth,trim={12cm, 0cm, 13cm, 0cm},clip]{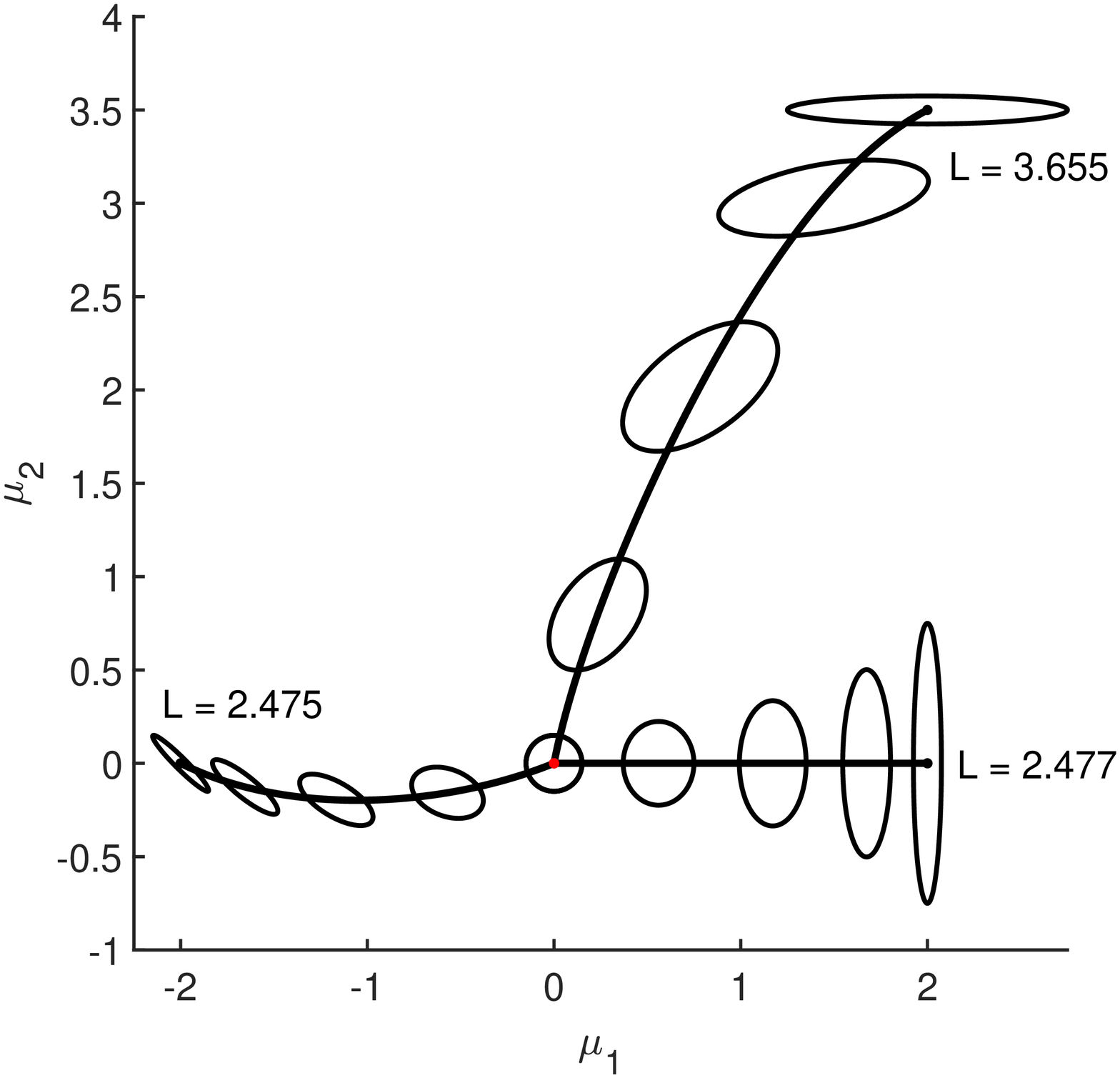}
    \includegraphics[width=0.48\textwidth,trim={13cm, 0cm, 13cm, 0cm},clip]{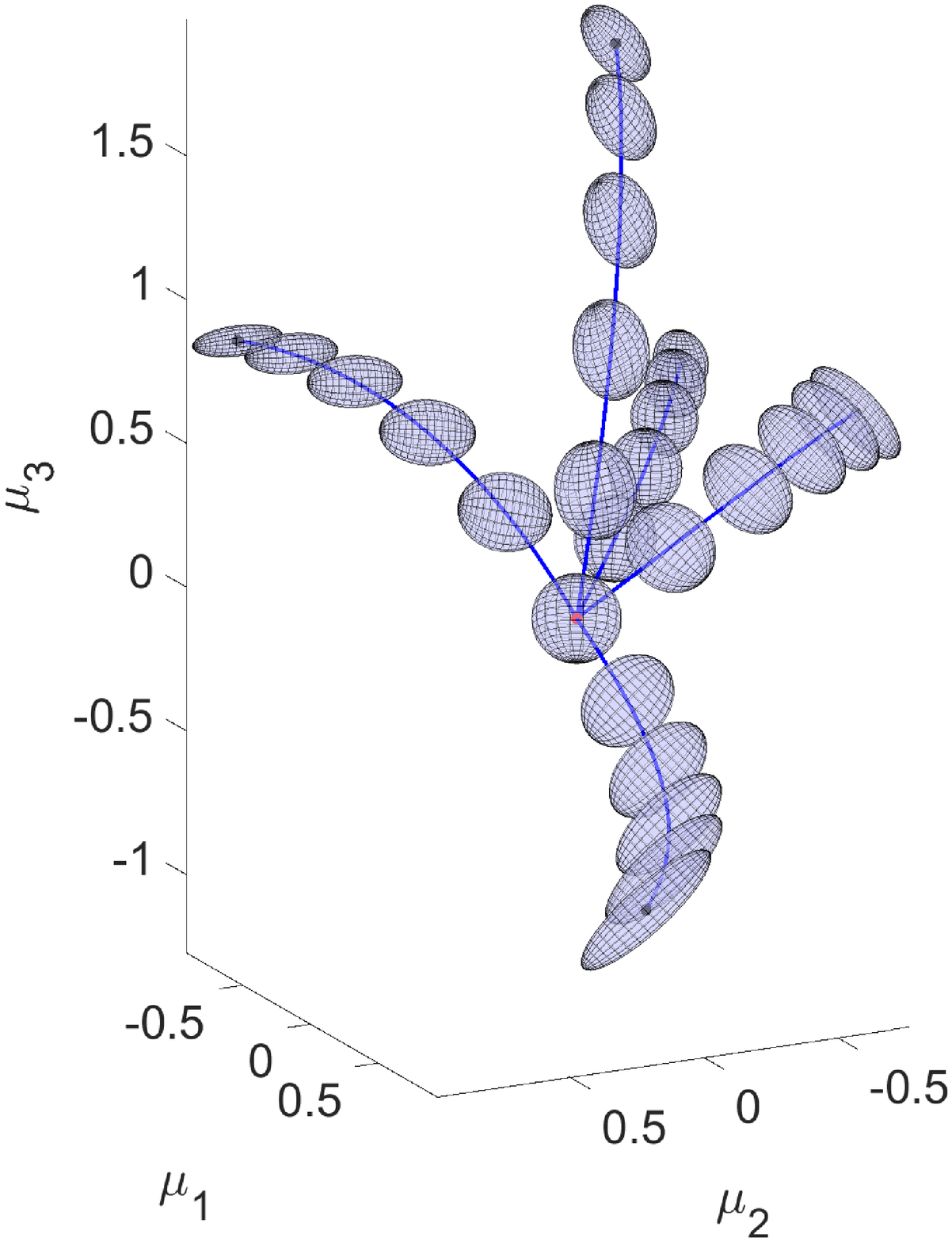}
    \caption{Example geodesics starting from the origin for bivariate (left) and trivariate (right) MVNs. Drawn paths indicate movement through $\bmu$-space, while positions in $\bSigma$-space at regular intervals along this path are indicated by (scaled) ellipses/ellipsoids. Geodesic paths curve through $\bmu$-space to optimise movement in higher-variance directions and often increase variance in the direction of travel.}
    \label{fig:geodesic_examples}
\end{figure}

Before describing the geodesics for the space of MVNs mathematically, we present in \figref{fig:geodesic_examples} some example geodesics in the space of bivariate and trivariate normal distributions that illustrate how the metric~\eqref{MVNmetric} influences their shapes. As discussed above, only the covariance matrix $\bSigma$ contributes to the warping of space, and so we may interpret $\bmu$-space as functioning in a Euclidean fashion outside of this ``external'' influence. As such, where we discuss ``straight'' or ``curved'' paths in $\bmu$-space we use these terms in the Euclidean sense (in the Riemannian sense, the geodesics function as the analog of straight lines by definition). Following others~\citep{Han2014, Pinele2020}, we also plot bivariate and trivariate MVNs and their geodesics in this sense, as points and paths through $\bmu$-space with covariance indicated by ellipses/ellipsoids with axis lengths given by (a scaling of) the eigenvalues of $\bSigma$, and orientations given by its eigenvectors.

Geodesic paths are seen to curve in $\bmu$-space in order to strike the balance between moving towards the target $\bmu$ value in the Euclidean sense, and travelling along higher-variance directions (which as per equation~\eqref{MVNmetric} incur less distance). In general, the geodesics initially increase the variance in the direction of travel, easing their movement through $\bmu$-space before restoring the covariance to that of the curve's endpoint. The endpoint's covariance can also bias the direction in $\bmu$-space from which the geodesic approaches, providing a motivation to leave the starting point from a direction that does not immediately point at the target. Notably however, when a principal direction of variance (as defined by the covariance matrix's eigenvectors) is aligned with the Euclidean connection between the two points in $\bmu$-space, the geodesic travels along this path. We return to this point in more detail in Section~\ref{sec:known_solutions}. 

Geodesics on Riemannian manifolds are found by solving the geodesic equations, which are ordinary differential equations (ODEs) defined wholly by the metric tensor. For the space of MVN distributions, these equations take the neat form~\citep{Imai2011}
\begin{equation}
\label{MVNgeodesics_base}
\begin{aligned}
\den{\bmu}{t}{2} &= \de{\bSigma}{t}\, \bSigma^{-1}\,  \de{\bmu}{t} \\
\den{\bSigma}{t}{2} &= \de{\bSigma}{t} \bSigma^{-1} \de{\bSigma}{t} - \left( \de{\bmu}{t} \right) \left( \de{\bmu}{t} \right)^T\!\!.
\end{aligned}
\end{equation}
In order to simplify the geodesic equations, it is beneficial to consider geodesics originating from the origin, $(\bmu, \bSigma) = (\b0, \bI)$. Making use of the freedom to transform according to~\eqref{affine_transform} without changing the metric, we may choose
\begin{equation}
    \label{affine_origin}
\bP = \bSigma^{1/2} \bR, \qquad \br = \bmu
\end{equation}
to map the origin, $(\b0, \bI)$, to any given point $(\bmu, \bSigma)$. Here, $\bR$ is an orthogonal matrix $\bR \bR^T = \bI$ that reflects the freedom to choose $\bP$ as any decomposition $\bP\bP^T = \bSigma$. This is equivalent to the freedom to choose from a variety of transformation matrices that all carry out the process of data whitening~\citep{Kessy2018}.

Considering \eqref{affine_origin}, we see that the geodesic between a starting point, $(\bmu_0, \bSigma_0)$, and a target point, $(\bmu_1, \bSigma_1)$, can instead be treated as a connection between the origin and a new transformed target,
\begin{equation}
\label{affine_target}
(\bmu_t, \bSigma_t) = \Bigl( \bR^T \bSigma_0^{-1/2} ( \bmu_1 - \bmu_0 ), \,\, \bR^T \bSigma_0^{-1/2} \bSigma_1 \bSigma_0^{-1/2} \bR \Bigr).
\end{equation}
The desired geodesic can then be obtained by applying the transform~\eqref{affine_origin} to this connection according to \eqref{affine_transform}. For a geodesic beginning at the origin, the geodesic equations~\eqref{MVNgeodesics_base} can be integrated once to become
\begin{equation}
\label{MVNgeodesic_odes}
\begin{aligned}
\de{\bmu}{t} &= \bSigma \bx          & \qquad \quad &  \bx = \left. \de{\bmu}{t} \right|_{t=0} \\
\de{\bSigma}{t} &= \bSigma( \bB - \bx \bmu^T )      & \qquad \quad  &  \bB = \left. \de{\bSigma}{t} \right|_{t = 0}.
\end{aligned}
\end{equation}
Several authors~\citep{Calvo1991,Imai2011,Eriksen1987} have shown that these equations permit a closed-form solution when expressed in terms of the canonical co-ordinates $(\bdelta, \bDelta) = (\bSigma^{-1}\bmu, \bSigma^{-1})$,
\begin{align}
\label{MVNgeodesic_solution}
    \bdelta(t) &= - \bB \Bigl( \cosh( t \bG) - \bI \Bigr) (\bG^{-})^{2}\bx + \sinh( t \bG ) \bG^{-} \bx, \\
    \nonumber
    \bDelta(t) &= \bI + \frac{1}{2} \Bigl( \cosh( t \bG ) - \bI  \Bigr) + \frac{1}{2} \bB \Bigl( \cosh(t \bG) - \bI \Bigr) (\bG^{-})^{2} \bB - \frac{1}{2} \sinh(t\bG) \bG^{-} \bB - \frac{1}{2} \bB \sinh(t\bG) \bG^{-}.
\end{align}
Here $\bG^2 = \bB^2 + 2\bx\bx^T$, and $\bG^{-}$ denotes the generalised inverse of $\bG$, satisfying $\bB \bG \bG^{-} = \bG \bG^{-} \bB = \bB$ and $\bG \bG^{-} \bx = \bG^{-} \bG \bx$ (see~\citep{Imai2011} for details). Notably, owing to the forms of the (globally convergent) Taylor series for the hyperbolic sine and cosine functions, it can be seen that equations \eqref{MVNgeodesic_solution} calls for only even powers of $\bG$, and so there is no degeneracy associated with knowing only $\bG^2$, and not $\bG$. 

Although equations~\eqref{MVNgeodesic_solution} can be evaluated directly, in our experience points along the geodesic path are best numerically evaluated using the form of the solution presented in~\citep{Eriksen1987}. In this form, $\bdelta(t)$ and $\bDelta(t)$ are found by reading off a portion of the result of exponentiating a matrix of augmented dimension,
\begin{equation}
    \label{MVNgeodesic_solution_matrix}
\begin{pmatrix}
\bDelta & \bdelta & \cdot\,\, \\
\bdelta^T & \cdot & \cdot\,\, \\
\cdot & \cdot & \cdot\,\,
\end{pmatrix} = 
\exp\left( t \begin{pmatrix} 
-\bB & \bx & \b0 \\
\bx^T & 0 & -\bx^T \\
\b0 & -\bx & \bB 
\end{pmatrix} \right).
\end{equation}
Here elements of the matrix marked with $\cdot$ are unused.

By observing that the velocity along geodesic curves is constant, the length of these curves over the segment $t \in [0,1]$ (the Fisher--Rao distance between the distributions at their start and endpoints) can be calculated by evaluating the velocity at the origin according to the metric~\eqref{MVNmetric}, giving~\citep{Pinele2020}
\begin{equation}
\label{MVNgeodesic_distance}
d_F = \sqrt{ \frac{1}{2} \tr(\bB^2) + \bx^T \bx } = \sqrt{ \frac{1}{2} \tr(\bG^2) }\,.
\end{equation}

Given a geodesic's initial velocity $\bv = (\bx,\bB)$, we thereby have closed-form expressions for both the path the geodesic takes~\eqref{MVNgeodesic_solution}, and the Fisher--Rao distance to the point at $t=1$~\eqref{MVNgeodesic_distance}. However, as yet there is no known analytic solution to the problem of choosing $\bv$ so as to produce a geodesic passing through a given arbitrary target point.  Outside of the special cases discussed subsequently, geodesics between arbitrary points in the MVN space can be found using a specially-tailored shooting method~\citep{Han2014}, as we discuss in Section~\ref{sec:shooting}. If only the Fisher--Rao distance between a pair of MVNs is required, closed-form bounds may suffice~\citep{Nielsen2023, Strapasson2016}.

\subsection{Special Cases with Known Solution}
\label{sec:known_solutions}

\subsubsection{Regions in $\bmu$-space Admitting Closed-Form Solution}
\label{sec:special_cases}
For several special cases, the velocity of the geodesic connecting a pair of MVNs is known, as summarised by~\citet{Pinele2020} and also~\citet{Nielsen2023} using the lens of totally geodesic submanifolds. Also relevant is the work of~\citet{Herntier2022}, who consider solutions for closest points on submanifolds (and geodesic paths to them), given some starting point on the full statistical manifold for bivariate MVNs. Here, we present the cases with known closed-form solution without any further invocation of differential geometry.

As equation~\eqref{MVNgeodesic_solution} gives the form of geodesics that begin at the origin, we discuss the cases with known solution in the context of having already been affine transformed so that they begin at the origin. That is, rather than the problem of finding the geodesic connecting two MVNs $(\bmu_0, \bSigma_0)$ and $(\bmu_1, \bSigma_1)$ we consider the analogous problem of finding the geodesic connecting the pair $(\b0, \bI)$ and $(\bmu_t, \bSigma_t)$, with the latter given by equation~\eqref{affine_target}. The desired initial velocity in this transformed space is $\bv = (\bx, \bB)$, and when it is obtained, the velocity of the geodesic connecting the actual targets is given by undoing the transform,
\[
\left.\dot{\bmu}\right|_{t=0} = \bSigma_0^{1/2} \bR \bx, \qquad \qquad \left.\dot{\bSigma}\right|_{t=0} = \bSigma_0^{1/2} \bR \bB \bR^T \bSigma_0^{1/2}.
\]
Owing to the invariance of the metric to these affine transforms, the Fisher--Rao distance between the original pair and the transformed pair remains the same.

We provide the cases with known solution in~\tabref{tab:known_geodesics}, and the associated Fisher--Rao distances in~\tabref{tab:known_distances}. Although these special cases may appear distinct, they are all in fact examples of a single type of solution. To demonstrate this, we observe that the choice of $\bR$ in equation~\eqref{affine_target} is free, and so we may choose $\bR$ so that it diagonalises $\bSigma_0^{-1/2}\bSigma_1\bSigma_0^{-1/2}$ and the target covariance becomes diagonal, $\bSigma_t = \bD$. Given this, applicability of the ``axis-aligned targets'' case in~Tables~\ref{tab:known_geodesics} and~\ref{tab:known_distances} depends only on its second condition, that the precision-weighted target mean, $\bdelta_t = \bSigma^{-1} \bmu_t$, is aligned with a principal axis after applying this choice of transformation,
\[
\bSigma_t^{-1} \bmu_t = \delta_t \be_{k}.
\]
Here $\delta_t = \|\bdelta_t\|_2$ is a scalar and $\be_{k}$ an elementary basis vector with $k$-th element unity and all other elements zero. This condition is satisfied only if $\bmu_t$ is aligned with an eigenvector of $\bSigma_t$.

\begin{table}[htp]
    \centering
    \begin{tabular}{@{}ccc@{}}
         \toprule \\[-0.5cm]
         {\bf Case} & $\bB$ & $\bx$ \\ \midrule \\[-0.4cm]
         Univariate, & \multirow{2}{*}{$\dfrac{\alpha}{\gamma}g$} & \multirow{2}{*}{\vspace{0.5cm}$\dfrac{2 g \delta_t \Delta_t}{\gamma}$} \\
         $\mu_t = \delta_t / \Delta_t,\,\, \Sigma_t = 1 / \Delta_t$ & & \\[0.25cm]
         Equal means, $\bmu_t = 0$ & $\log \left( \bSigma_t \right)$ & $\b0$ \\[0.4cm]
         Equal variances, $\bSigma_t = \bI$ & $\dfrac{\acosh\left(1 + \|\bmu_t\|^2 + \frac{1}{8}\|\bmu_t\|^4\right)}{\|\bmu_t\| \sqrt{\|\bmu_t\|^2 + 8}}\,\, \bmu_t \bmu_t^T$ & $\dfrac{2 \acosh\left(1 + \|\bmu_t\|^2 + \frac{1}{8}\|\bmu_t\|^4\right)}{\|\bmu_t\| \sqrt{\|\bmu_t\|^2 + 8}} \,\, \bmu_t$ \\[0.6cm]
         Axis-aligned targets, & \multirow{2}{*}{$\Bigl( \dfrac{\alpha}{\gamma} g\, \bI - \log(\bD) \Bigr) \be_k \be_k^{T} + \log(\bD)$} & \multirow{2}{*}{\vspace{0.25cm}$ \dfrac{2g\delta_t \Delta_t}{\gamma}\,\, \be_k $} \\
         $\bSigma_t^{-1}\bmu_t = \delta_t \be_{k},\,\, \bSigma_t = \bD$ & & \\[0.2cm] \midrule
         \multicolumn{3}{c}{Constants: \quad $\Delta_t = 1/D_{kk},\quad \alpha = \delta_t^2 - 2 \Delta_t^2 + 2 \Delta_t, \quad \gamma = \sqrt{\alpha^2 + 8 \delta_t^2 \Delta_t^2} $, \quad $g = \acosh \left(1 + \gamma^2 / 8 \Delta_t^3 \right)$} \\
         \bottomrule
    \end{tabular}
    \caption{Special cases for which the geodesic connecting the origin to a target point, $(\bmu_t,\bSigma_t)$, is known. The target point $(\bmu_t, \bSigma_t)$ associated with a given start and end point is provided by equation~\eqref{affine_target}. The condition that the target mean is aligned with an elementary basis vector, $\bSigma^{-1}_t \bmu_t = \delta_t \be_k$ defines the index $k$, and for the univariate case $k = 1$. $\bD$ denotes a diagonal matrix.}
    \label{tab:known_geodesics}
\end{table}

\begin{table}[htp]
    \centering
    \begin{tabular}{@{}cc@{}}
         \toprule \\[-0.5cm]
         {\bf Case} & $d_F$ \\ \midrule & \\[-0.4cm]
         Univariate, & \multirow{2}{*}{$\dfrac{1}{\sqrt{2}} \acosh\left(1 + \dfrac{(\delta_t^2 - 2\Delta_t^2 + 2\Delta_t)^2 + 8\delta_t^2\Delta_t^2}{8 \Delta_t^3} \right)$} \\
         $\mu_t = \delta_t / \Delta_t,\,\, \Sigma_t = 1 / \Delta_t$ & \\[0.4cm]
         Equal means, $\bmu_t = 0$ & $\sqrt{\dfrac{1}{2}\sum_{i=1} (\log \lambda_i)^2}, \qquad \mbox{$\lambda_i$ the eigenvalues of $\bSigma_t$}$ \\[0.5cm]
         Equal variances, $\bSigma_t = \bI$ & $\dfrac{1}{\sqrt{2}} \acosh\left( 1 + \|\bmu_t\|^2 + \frac{1}{8} \|\bmu_t\|^4 \right)$ \\[0.6cm] 
         Axis-aligned targets, & \multirow{2}{*}{$\sqrt{ \frac{1}{2} \acosh^2\left(1 + \dfrac{(\delta_t^2 - 2D_{kk}^2 + 2D_{kk})^2 + 8\delta_t^2D_{kk}^2}{8 D_{kk}^3} \right) + \frac{1}{2}\sum\limits_{j \neq k} (\log D_{jj})^2}$} \\[0cm]
         $\bSigma_t^{-1}\bmu_t = \delta_t \be_{k},\,\, \bSigma_t = \bD$ & \\[0.35cm] 
         \bottomrule
    \end{tabular}
    \caption{Fisher--Rao distances $d_F$ associated with the special cases given in \tabref{tab:known_geodesics}. The target point $(\bmu_t, \bSigma_t)$ associated with a given start and end point is provided by equation~\eqref{affine_target}. The condition that the target mean is aligned with an elementary basis vector, $\bSigma^{-1}_t \bmu_t = \delta_t \be_k$ defines the index $k$, and for the univariate case $k = 1$. $\bD$ denotes a diagonal matrix.}
   \label{tab:known_distances}
\end{table}

Such a condition is satisfied by all of the simpler cases in Tables~\ref{tab:known_geodesics} and~\ref{tab:known_distances}. For the univariate case, all 1-vectors are aligned by definition and the condition is trivially satisfied. For the case of equal means, $\bmu_t = 0$ and the zero-length vector can be considered aligned with any eigenvector of the target covariance. For the case of equal covariance matrices, $\bSigma_t = \bI$ and choice of $\bR$ is not restricted to be the diagonaliser. As such, we can instead choose $\bR$ that rotates $\bmu_t$ to align with a co-ordinate axis and hence satisfy the condition.

\begin{figure}
    \centering
    \includegraphics[width=0.75\textwidth,trim={18cm, 8cm, 16cm, 6.5cm},clip]{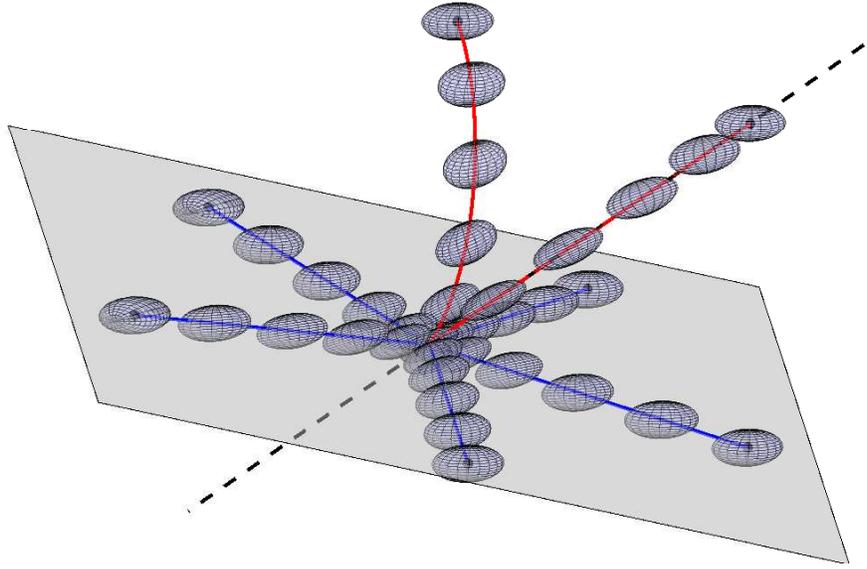}
    \caption{Spaces of target means $\mu_1$ reachable by ``straight line'' geodesics, as imposed by the eigenstructure of $\bSigma_t = \bSigma_0^{-1/2} \bSigma_1 \bSigma_0^{-1/2}$. Pictured are a set of example geodesics for trivariate MVNs, each starting from the same point $(\bmu_0, \bSigma_0)$ and sharing the same $\bSigma_1$. Here, $\bSigma_t$ has a repeated eigenvalue (multiplicity 2) and produces a plane of $\bmu_1$ values for which geodesic paths (blue) travel straight lines through $\bmu$ space with closed form (\tabref{tab:known_geodesics}). Geodesics to points with $\bmu_1$ outside of the plane (red) only take this form if they align with the remaining eigenvector of $\bSigma_t$ (dotted line).}
    \label{fig:eigenspace_examples}
\end{figure}

More precisely, the requirement for the use of this solution is that $\bmu_t$ lies within the eigenspace associated with one of the eigenvalues of $\bSigma_t$. Given that the geometric and algebraic multiplicities of all covariance matrix eigenvalues are equal, each eigenvalue repeated $n$ times thus corresponds to a (hyper)plane of dimension $n$ in which $\bmu_t$ (after rotating to diagonalise $\bSigma_t$) may lie for a closed-form solution to be available. \figref{fig:eigenspace_examples} presents a graphical demonstration of this concept, showing how in the space of trivariate MVNs, an eigenvalue of multiplicity two and another of multiplicity one correspond to a plane and a line, respectively, of target means $\bmu_t$ permitting closed-form geodesics. 

The closed-form geodesics travel straight lines in $\bmu$ space. Mathematically, these correspond to the cases where $\bB$ and $\bG$ are mutually diagonalisable, allowing for each dimension in equations~\eqref{MVNgeodesic_solution} to be decoupled via a similarity transform (see supplementary material). In these cases, the covariance matrix eigenvectors remain fixed in specific directions along the geodesic, and geodesics travel along straight lines oriented in one of these directions. This property is key to a geometric understanding of these solutions, and helps understand the geodesics of the MVN statistical manifold more generally.

As the metric~\eqref{MVNmetric} dictates, increasing the variance in a direction lowers the information cost of moving in that direction. Aligning the covariance eigenvectors with the direction of movement thus becomes a way to potentially decrease overall information cost, by providing the ability to increase variance only in the direction where this will be useful. General geodesics must balance this effect against the orientation of eigenvectors at the geodesic's start and end, as well as the curved path they take through $\bmu$-space. However, when the direction to the target already aligns with a covariance eigenvector, there is no longer any motivation to rotate the principal directions of variance. The geodesic solution optimises for increasing the variance associated with that direction to ease movement, and reaching the correct target, while all other eigenvalues have no movement through $\bmu$-space to consider and simply transform according to the equal means solution, such that $\log \lambda$ linearly approaches the target covariance eigenvalue $\lambda_t$,
\begin{equation}
\label{exponential_variance}
\lambda(t) = (\lambda_t)^t.
\end{equation}

\begin{figure}
    \centering
    \includegraphics[width=0.5\textwidth,trim={0.5cm, 0cm, 2cm, 1cm},clip]{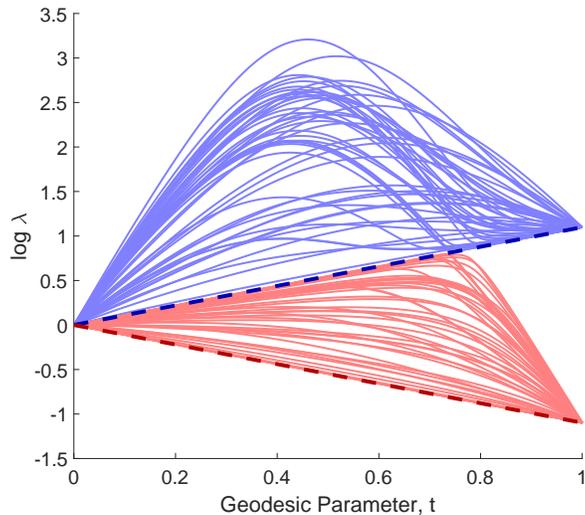}
    \caption{Paths traced out by the larger (blue) and smaller (red) eigenvalue of $\bSigma$, for geodesics in the space of bivariate MVNs connecting the origin to random points with final covariance eigenvalues $\log \lambda_1 = 1$, $\log \lambda_2 = -1$. The straight dotted lines depict $\log \lambda = (\log \lambda_t) t$, the lower bound for eigenvalues along all geodesic paths.}
    \label{fig:eigenvalue_bound}
\end{figure}

Even when the direction to the target in $\bmu$-space does not align with an eigenvector of the target covariance, equation~\eqref{exponential_variance} still serves as a lower bound for the values of covariance eigenvalues even while their associated eigenvectors shift direction (\figref{fig:eigenvalue_bound}). This is geometrically understood by considering that equation~\eqref{exponential_variance} represents the lowest-cost path eigenvalues follow to reach a given target, ignoring influence from any movement through $\bmu$-space. When we do consider that movement, there may be a motivation to increase variance in that direction in order to lower the information cost associated with it. However, there is never any reason to {\it decrease} variance below the path defined by equation~\eqref{exponential_variance}. Interestingly, this result emerges naturally from our geometric understanding, but seems far more difficult to demonstrate mathematically.

\section{Non-Geodesic Paths}

\subsection{Commonly-Used Paths}
\label{sec:common_paths}

Given the observation that paths with low information cost improve samplers using distribution sequences~\citep{Sim2012, Syed2021}, we summarise here how commonly-used paths in such settings function in terms of information cost, using our distribution family of MVNs as an accessible way to understand how these different choices of paths behave. Some of these paths along the MVN manifold were discussed in a very recent publication~\citep{Nielsen2023}, however in keeping with the character of this review we concentrate on practical performance and an interpretation external to information geometry. For notational consistency with the geodesic paths already discussed, we parameterise these paths using $t \in [0,1]$, with $t = 0$ corresponding to the starting point $\btheta_0 = (\bmu_0, \bSigma_0)$ and $t = 1$ to the endpoint, $\btheta_1 = (\bmu_1, \bSigma_1)$.

A predominant choice in a variety of contexts is the annealing path~(e.g. \citep{Hukushima1996,Duan2015}), which takes the form
\[
p(t;\bx) = \bigl[p_0(\bx)\bigr]^{1-t} \,\bigl[p_1(\bx)\bigr]^t, \qquad \quad \mbox{equivalently} \quad \begin{aligned} \bdelta(t) &= (1 - t) \bdelta_0 + t \bdelta_1 \\ \bDelta(t) &= (1 - t) \bDelta_0 + t \bDelta_1. \end{aligned}
\]
This path represents a linear interpolation between the log densities associated with the two distributions, and as MVNs are a member of the exponential family, the annealing path also corresponds to a linear interpolation in terms of their natural parameters~\citep{Grosse2013}, and equivalently, $\bdelta$ and $\bDelta$. Observing that the geodesic paths~\eqref{MVNgeodesic_solution} are non-linear in these canonical co-ordinates immediately demonstrates that annealing paths are non-geodesic.

A second choice of path for distributions in a single exponential family is the moment-averaged path, which can offer significantly improved performance in annealed importance sampling~\citep{Grosse2013}. For MVNs, this path takes the form
\begin{equation}
\label{moment_path}
\begin{aligned}
\bmu(t) &= (1 - t) \bmu_0 + t \bmu_1 \\
\bSigma(t) &= (1 - t) \bSigma_0 + t \bSigma_1 + \frac{1}{4} t (1 - t) (\bmu_1 - \bmu_0) (\bmu_1 - \bmu_0)^T.
\end{aligned}
\end{equation}
These paths move in Euclidean fashion through $\bmu$-space, which immediately renders them non-geodesic. However, these paths hold the attractive property that this movement in $\bmu$-space is eased by the temporary injection of variance in that direction, mimicking one of the key behaviours of geodesic paths.

A somewhat separate type of path are those defined by the Wasserstein distance. The Wasserstein distance is formally a distance with an associated metric~\citep{Panaretos2019}, albeit one of different form to that of the Fisher metric. Informally, the Wasserstein distance measures the amount of physical work required to move the probability mass of one density function, in the fashion that minimises the work performed, so that it becomes the density function of another. Paths that are geodesic with respect to the Wasserstein metric are not in general geodesic with respect to the Fisher metric, but are interesting to consider given the use of Wasserstein distance as a statistical distance measure~\citep{Arjovsky2017} and in averaging distributions via the Fr{\'{e}}chet mean~\citep{Cuturi2014}.

In the space of MVNs, the 2-Wasserstein distance is typically used as it provides a closed form for the geodesic path (and distance) between two MVNs~\citep{Delon2020}. This path takes the form
\begin{equation}
\label{wasserstein_path}
\begin{aligned}
\bmu(t) &= (1 - t) \bmu_0 + t \bmu_1 \\
\bSigma(t) &= \Bigl( (1 - t) \bI + t \bC \Bigr) \bSigma_0 \Bigl( (1 - t) \bI + t \bC \Bigr), \qquad \bC = \bSigma_1^{1/2} \left( \bSigma_1^{1/2} \bSigma_0 \bSigma_1^{1/2} \right)^{-1/2} \bSigma_1^{1/2}.
\end{aligned}
\end{equation}
Paths minimising the 2-Wasserstein distance also move in Euclidean fashion through $\bmu$-space, but without any explicit injection of variance --- equations~\eqref{wasserstein_path} may be interpreted as applying a linearly-varying affine transformation to the starting point of the path, starting with the identity transformation and ending with an affine transformation that maps $(\bmu_0, \bSigma_0)$ to $(\bmu_1, \bSigma_1)$. As such, the Wasserstein metric does not directly acknowledge the information about the random variable contained in the distribution~\citep{Khan2022}, and indeed for the univariate case equations~\eqref{wasserstein_path} reduce to simply a linear interpolation of mean and standard deviation, $\sqrt{\Sigma(t)} = (1-t) \sqrt{\Sigma_0} + t \sqrt{\Sigma_1}$~\citep{Li2023}. Despite this information-unaware property, the Wasserstein distance has numerous applications in statistics, many of which are reviewed by~\citet{Panaretos2019}. It is also well suited to comparing distributions with disjoint support, in contrast to information theoretic approaches~\citep{Khan2022}.

\begin{figure}
    \hspace{0.025\textwidth}
    % FIRST FIGURE
    \begin{minipage}{0.025\textwidth}
    {\bf a)}
    \end{minipage}%
    \begin{minipage}{0.425\textwidth}
    \centering
    \includegraphics[width=0.95\textwidth, trim={13cm, 0cm, 11cm, 0cm}, clip]{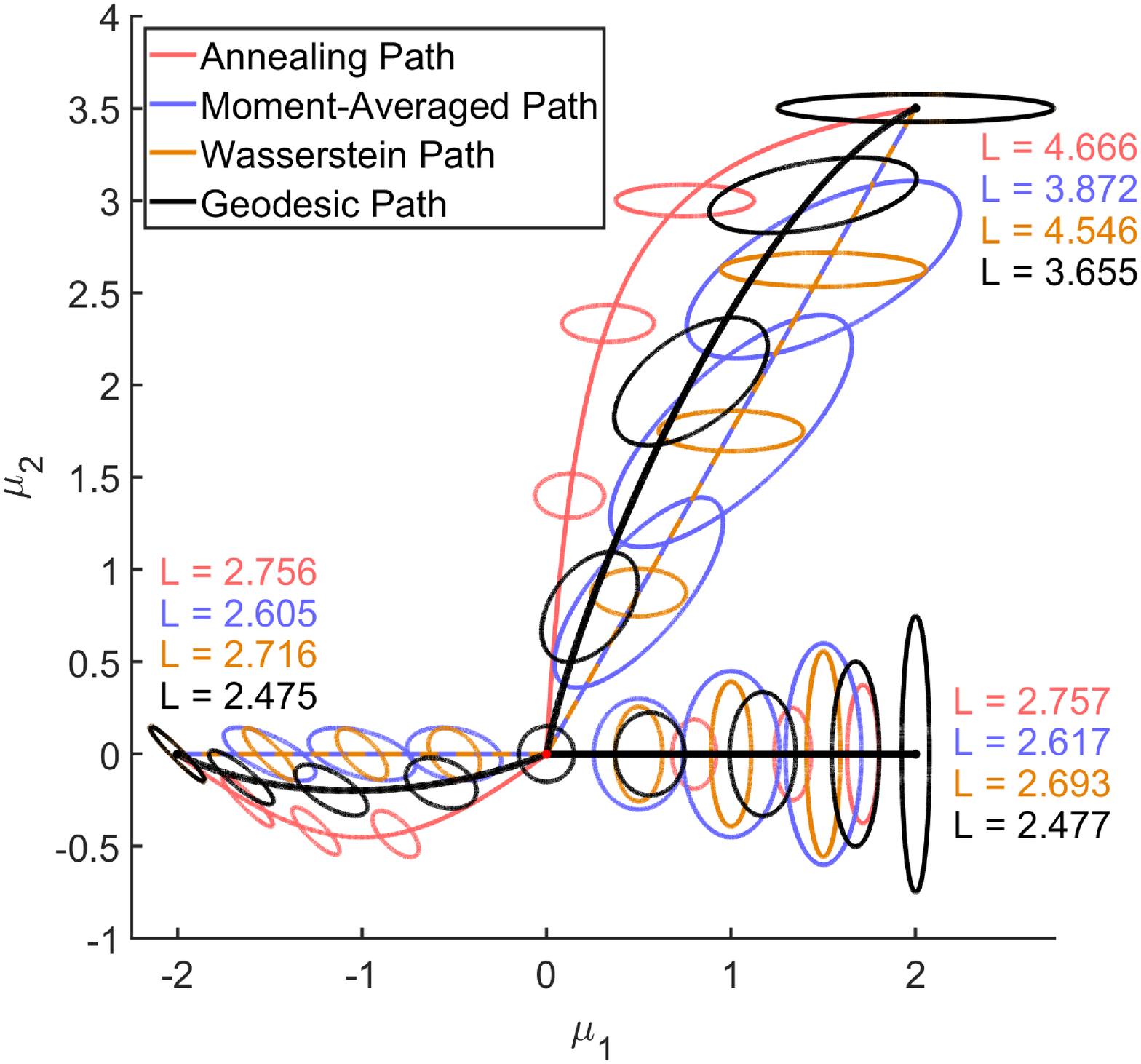}
    \end{minipage}%
    \hspace{0.04\textwidth}
    % SECOND FIGURE
    \begin{minipage}{0.025\textwidth}
    {\bf b)}
    \end{minipage}%
    \begin{minipage}{0.425\textwidth}
    \centering
    \includegraphics[width=0.95\textwidth, trim={1cm, 0cm, 3.5cm, 0cm}, clip]{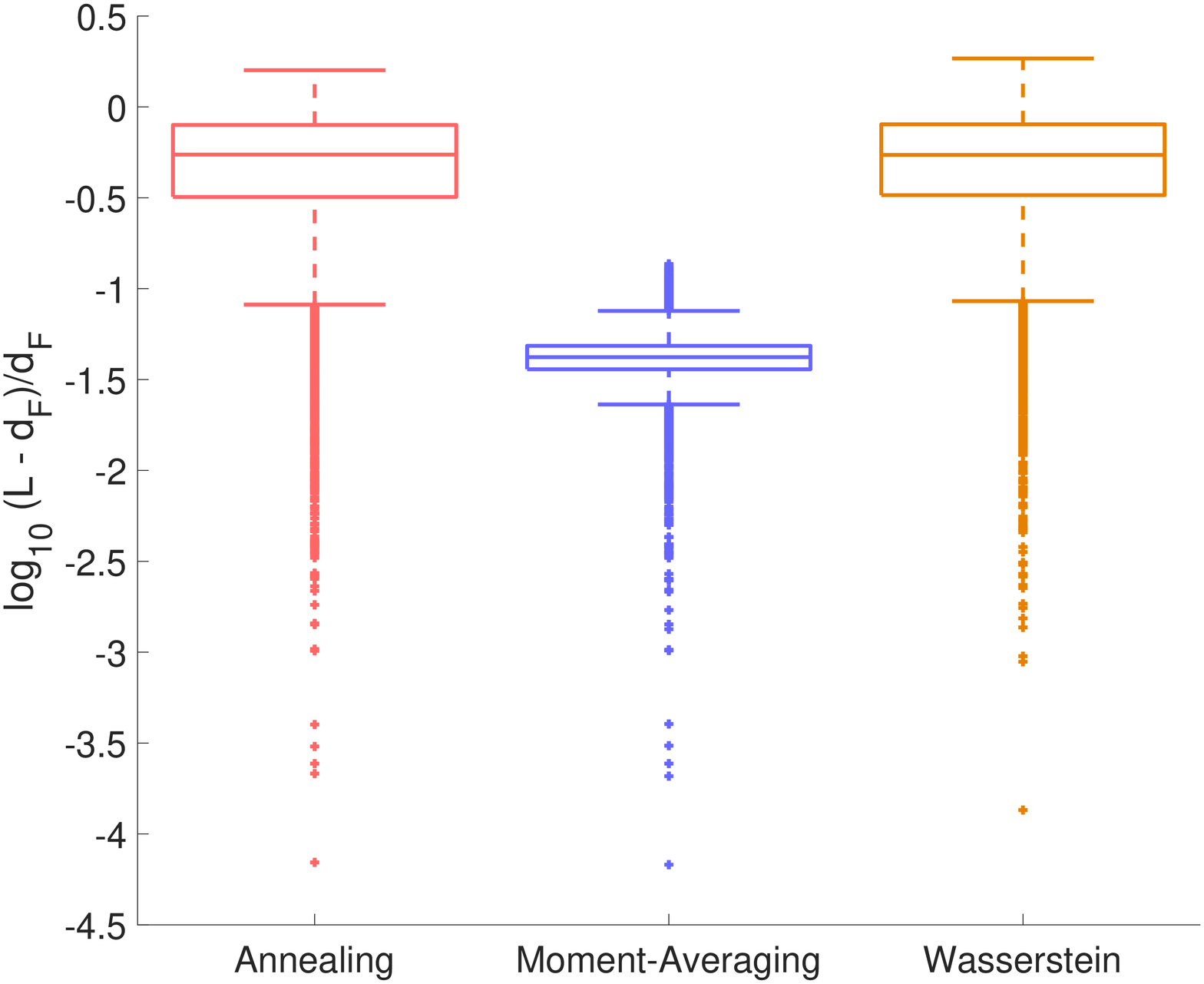}
    \end{minipage}%
    \caption{ {\bf a)} The paths taken through the space of bivariate MVNs, and associated manifold distances, of commonly-used non-geodesic paths for the examples shown in~\figref{fig:geodesic_examples}. %Annealing paths respond too strongly to the initial and final covariance matrices, resulting in overly curved paths through $\bmu$-space. In contrast, moment-averaging paths move in direct Euclidean fashion through $\bmu$-space, but increase variance in this direction to reduce the distance incurred by this choice.
    {\bf b)} Performance (in terms of manifold distance relative to the geodesic distance $d_F$) for different choices of path connecting the origin to a large number of random points in the space of MVNs. %Annealing and 2-Wasserstein-geodesic paths perform similarly across these tests despite taking very different trajectories. Moment-averaging paths of form~\eqref{moment_path} easily outperform both, with a median relative extra length of 4.2\% as compared to 54.6\% (annealing) and 54.4\% (Wasserstein).
    Moment-averaged paths prove superior.
    }
    \label{fig:MVN_paths}
\end{figure}

In~\figref{fig:MVN_paths}a, we visualise these different choices of closed-form path for a few examples in the space of bivariate MVNs. We immediately see that 2-Wasserstein-optimal paths move in direct fashion through the space, and are perhaps best interpreted as optimal {\it absent} the informational context. Annealing paths strongly favour travelling the directions of high variance at the start and end point, as they do not inject additional variance. These paths curve too much in $\bmu$-space relative to the geodesic path and hence become information-inefficient. Conversely, moment-averaged paths inject additional variance to ease movement through $\bmu$-space, but do not curve this movement to take advantage of the variance at the path's endpoints.

As observed by \citet{Grosse2013}, the behaviours of annealing and moment-averaged paths can be understood through a variational lens, as the former can be shown to be reverse-KL optimal (mode-seeking, avoiding additional variance) and the latter forwards-KL optimal (mean-seeking, adding additional variance). \citet{Nielsen2023} demonstrates the potential of paths formed by combining these two paths, trading off their respective strengths and weaknesses.

To test each type of path's performance more generally, we generate a large number of random points in MVN space by independently sampling components of $\bmu$ uniformly, $\mu_i \sim U[-10,10]$, sampling covariance eigenvalues log-uniformly, $\log \lambda \sim U[-10, 10]$, and creating covariance matrices from these eigenvalues by generating random orthonormal matrices through the Gram--Schmidt process applied to randomly selected direction vectors. We connect each of these points to the origin with the different types of paths, and calculate their length by summing the lengths of each of their segments~\eqref{MVNmetric} after applying a fine discretisation. Boxplots of these lengths across the many test problems are given in~\figref{fig:MVN_paths}b. Despite taking very differently-shaped paths through the space of MVNs, annealing paths and 2-Wasserstein-optimal paths perform extremely similarly in terms of the information distance associated with their traversal, at least for test problems generated in this manner. Moment-averaged paths are much better at travelling close-to-information-optimal paths through MVN space, in accordance with their demonstrated ability to improve annealed importance sampling~\citep{Grosse2013}.

\subsection{MVN-Specific Approximations}
\label{sec:approx_solutions}

Ideas that approximate the geodesic path, either by approximating the correct choice of $\bv$ to reach a given target or by generating paths that reach the target and are close to geodesic, are useful both for approximating or bounding the Fisher--Rao distance~\citep{Strapasson2016,Nielsen2023}, and for improving numerical efforts to find geodesic paths~(Section~\ref{sec:shooting}). Here we briefly discuss both types of approximation, in all cases considering the case where equation~\eqref{affine_target} has already been used to convert the problem of connecting two MVNs to the equivalent problem of connecting the origin to a given target, $\btheta_t = (\bmu_t, \bSigma_t)$.

\citet{Nielsen2023} observes that approximately-geodesic paths can be found by projecting from a higher-dimensional Riemannian manifold for symmetric positive definite matrices, $\cP$, for which \citet{Calvo1991} demonstrated a submanifold is isometric to the MVN Riemannian manifold. Geodesics connecting the origin to where the target point is embedded in $\cP$ are given simply by
\begin{equation}
\label{P_geodesic}
\bp(t) = \begin{bmatrix} \bSigma_t + \bmu_t \bmu_t^T & \bmu_t \\ \bmu_t^T & 1 \end{bmatrix}^{t} = \begin{bmatrix} p_{\bSigma}(t) & p_{\bmu}(t) \\ p_{\bmu}(t)^T & p_{\beta}(t) \end{bmatrix},
\end{equation}
where the latter matrix is defined purely for notational convenience. The sub-manifold corresponding to MVNs requires $p_{\beta} = 1$, but $\bp(t)$ violates this for intermediate values of $0 < t < 1$ and so leaves the MVN manifold. Nevertheless, approximately-geodesic paths through MVN space can be obtained by projecting back onto the MVN manifold~\citep{Nielsen2023}, via
\begin{equation}
    \label{projection_path}
\bmu(t) = \frac{p_{\bmu}(t)}{p_{\beta}(t)}, \qquad \bSigma(t) = \bp_{\bSigma}(t) - \frac{p_{\bmu}(t) p_{\bmu}(t)^T}{\bp_{\beta}(t)}.
\end{equation}
These paths are generally highly competitive, although \citet{Nielsen2023} observed in numerical experiments that for higher dimensions ($d > 11$) moment-averaged paths~\eqref{moment_path} became superior. We have also observed this to be the case for geodesics between points separated very far in (variance-weighted) $\bmu$-space, even in lower dimensions.

We may also approximate geodesics in the sense of finding choices for $\bv = (\bx, \bB)$ that approximate those of the desired geodesic, and these prove useful for numerically finding the true geodesic. Some of these approximations can be obtained by extending the special cases in~\tabref{tab:known_geodesics} beyond their range of application. The ``equal means'' solution corresponds to the case $\bB = \bG$ and as detailed in the supplement, we gain an approximation for small $\bx$,
\begin{equation}
    \label{taylor_approximation}
\begin{aligned}
\bx &\approx \Bigl( \bI - (\bI - \bSigma_t)^{-1} \Bigr) \log( \bSigma_t) \bSigma^{-1}_t \bmu_t\\
\bB &\approx \log \bSigma_t,
\end{aligned}
\end{equation}
by substituting $\bB^2 \approx \bG^2$ into equations~\eqref{MVNgeodesic_solution}.

We also obtain an approximate solution by relaxing the requirement that $\bmu_t$ be aligned with a covariance matrix eigenvector. Rather, we project the target mean onto each eigenvector separately, form the optimal univariate solution for each component and combine them,
\begin{equation}
    \label{MVNgeodesics_approx}
\bB \approx \sum_{j=1}^d \frac{\alpha_j}{\gamma_j} g_j \bv_j \bv_j^T, \qquad \quad \bx \approx \sum_{j=1}^d \sign(\bv_j^T \bdelta_t) g_j \sqrt{\frac{1}{2}\left(1 - \left(\frac{\alpha_j}{\gamma_j}\right)^2\right)} \bv_j,
\end{equation}
where $(\lambda_j, \bv_j)$ are the eigenpairs of $\bSigma_t^{-1}$, $\bdelta_t = \bSigma_t^{-1} \bmu_t$ as before, and
\[
\alpha_j = (\bv_j^T \delta)^2 - 2 \lambda_j^2 + 2 \lambda_j, \quad \gamma_j = \sqrt{ \alpha_j^2 + 8 (\bv_j^T \bdelta_t) \lambda_j^2 }, \quad g_j = \acosh \left( 1 + \gamma_j^2 / 8 \lambda_j^3 \right).
\]
This approximation is exact for any of the cases in~\tabref{tab:known_geodesics}. A related idea was used to bound the Fisher--Rao distance by \citet{Strapasson2016}.

Finally, we may also generate approximate initial velocities by considering the initial velocity of the geodesics~\eqref{P_geodesic} on $\cP$. As we demonstrate in the supplementary material, some blocks of the logarithm of the matrix in equation~\eqref{P_geodesic} provide an approximation for the required $\bx$ and $\bB$,
\begin{equation}
    \label{projection_approximation}
\bp'(0) = \log \begin{bmatrix} \bSigma_t + \bmu_t \bmu_t^T & \bmu_t \\ \bmu_t^T & 1 \end{bmatrix} \approx \begin{bmatrix} \bB & \bx \\ \bx^T & \cdot \end{bmatrix}.
\end{equation}
As with the small-$\bx$ approximation above, this approximation is exact in the case of equal means, as can be seen by direct substitution of $\bx = \b0$ and $\bmu_t = \b0$.

Our focus is on the use of these approximations inside the shooting algorithms discussed subsequently. However, in the supplementary we display a few example geodesics and the corresponding approximate geodesics (\figref{suppfig:approximate_geodesic_examples}) and their general performance across randomly-generated test problems~(\figref{suppfig:approximate_geodesic_performance}). The projection-based approximation~\eqref{projection_approximation} emerges as the best overall performer, especially in higher dimensions.

\section{Shooting to Find Geodesic Paths}
\label{sec:shooting}

\subsection{Fundamental Concepts}

Outside the cases listed in \tabref{tab:known_geodesics}, the choice of $\bx$ and $\bB$ in equation~\eqref{MVNgeodesic_solution} that produce the geodesic that reaches a given target point remains unknown. Instead, geodesics between arbitrary points in MVN space are produced numerically, via shooting~\citep{Han2014}. Shooting is the process of taking an initial guess for the velocity $\bv = (\bx, \bB)$, firing out the geodesic with this velocity, and then adjusting this velocity according to where the geodesic ray ends up relative to the target point. These adjustments are made leveraging the known geometric properties of the space. Convergence of shooting is nontrivial, and \citet{Han2014} suggest that for manifold distances greater than about seven the algorithm demonstrates issues with both stability and convergence rate. For these more challenging problems, they propose an alternative shooting via path refinement (Section~\ref{sec:path_refinement}).

Here, we first review the fundamental pieces of the algorithm of \citet{Han2014} to complete this work's differential geometry-free treatment of these statistical manifolds. Then, in Section~\ref{sec:improved_shooting}, we demonstrate the use of MVN-specific approximations (Section~\ref{sec:approx_solutions}) inside the shooting process.

After shifting a problem to the origin using equation~\eqref{affine_origin}, we can obtain the geodesic associated with a given initial velocity using equation~\eqref{MVNgeodesic_solution_matrix}. To formulate some correction to this velocity, we require a ``residual'' vector that compares the geodesic's endpoint and the target point. This vector becomes an update to the geodesic's initial velocity via the differential geometric process known as parallel transport, which moves a vector along a Riemannian manifold whilst transforming it so as to keep it pointing in the same direction relative to the manifold.

The parallel transport equations for the case of MVNs are given in~\citet{Han2014}, and for geodesics beginning at the origin they take the form
\begin{equation}
    \label{parallel_transport}
\begin{aligned}
\de{\bumu}{t} &= \frac{1}{2} \bSigma \left( \bB - \bx \bmu^T \right) \bSigma^{-1} \bumu + \frac{1}{2} \buSigma \, \bx \\
\de{\buSigma}{t} &= \frac{1}{2} \bSigma \left( \bB - \bx \bmu^T \right) \bSigma^{-1} \buSigma + \frac{1}{2} \buSigma \left(\bB - \bx \bmu^T \right) - \frac{1}{2} \bSigma \, \bx\,  \bumu^T - \frac{1}{2} \bumu \, \bx^T \, \bSigma.
\end{aligned}
\end{equation}
Integrating equations~\eqref{parallel_transport} backwards along a geodesic (that is, from $t = 1$ to $t = 0$), we can thus shift the residual vector back to the start of the geodesic, where it can be used to update our initial velocity guess. Integration of these equations is achieved numerically, with the values of $\bmu$ and $\bSigma$ at each timepoint given by~\eqref{MVNgeodesic_solution}.

The issue is, for a non-Euclidean manifold the traditional concept of a residual vector is unavailable. \citet{Han2014} use the Euclidean residual as a ``good enough'' choice in their shooting approach, reducing the magnitude of resulting velocity updates when the magnitude of this vector grows too large and the resulting approximation unsatisfactory. Observing that the residual vector in Euclidean space is also the velocity of the geodesic connecting the two points (a straight line), here we also consider the approximate velocities in Section~\ref{sec:approx_solutions} as alternative choices of residual vector (Section~\ref{sec:improved_shooting}). Note, however, that even if the exact geodesic velocity is available, this does not produce a perfect update.

In order to decide an appropriate step length for a velocity update, shooting makes use of the Riemannian geometry concept known as the Jacobi field. The Jacobi field describes the difference between infinitesimally close geodesics, and hence in a fashion analogous to the derivative or Jacobian in other contexts, can be used to obtain the linearised sensitivity of a geodesic to some perturbation. To explicitly define the Jacobi field and the numerical shooting algorithms, we first introduce the compact notation $\cG(t \bv)$ to denote the action of tracing out a geodesic through MVN space from the origin $(\b0, \bI)$ with initial velocity $\bv = (\bx, \bB)$ out to distance $t$. This function $\cG(t\bv)$ is given by either of equations~\eqref{MVNgeodesic_solution} or~\eqref{MVNgeodesic_solution_matrix}. The problem of finding the geodesic between two points in MVN space may thus be expressed as seeking an initial velocity $\bv$ such that
\[
\btheta_t = \cG(\bv), \qquad \qquad \mbox{or equivalently,} \quad \bv = \cG^{-1}(\btheta_t),
\]
where again $\btheta_t = (\bmu_t, \bSigma_t)$ is the target point after shifting to the origin. We note that in differential geometry, $\cG(\bv)$ and $\cG^{-1}(\btheta)$ are often described as the exponential map, $\exp_{(\b0,\bI)}(\bv)$, and the logarithmic map, $\log_{(\b0,\bI)}(\btheta)$, respectively.

Using this notation, the Jacobi field associated with a given velocity update $\Delta \bv$ is given by
\begin{equation}
    \label{jacobi_field}
    \cJ (t) = \left.\deb{\tau} \cG\Bigl(t (\bv + \tau \Delta \bv) \Bigr)\right|_{\tau = 0} = \lim_{\epsilon \rightarrow 0} \frac{\cG\Bigl(t (\bv + \epsilon \Delta \bv)\Bigr) - \cG(t\bv)}{\epsilon}.
\end{equation}
Estimating $\cJ (1)$ using finite differencing, we obtain the linearised effect of a velocity update, $\Delta \bv$, on the end point of the geodesic,
\begin{equation}
\label{Jacobi1}
\cJ (1) \approx \frac{\cG(\bv + s \Delta \bv) - \cG(\bv)}{s}.
\end{equation}
Seeking a scalar stepsize $s$ such that $\cG(\bv + s \Delta v) \approx \btheta_t$, we can substitute this relation into equation~\eqref{Jacobi1} and then scalarise it by taking an inner product of both sides,
\begin{equation}
    \label{stepsize}
s \approx \frac{\Bigl\langle \btheta_t - \cG(\bv), \cJ (1) \Bigr\rangle_{\cG(\bv)}}{\Bigl\langle \cJ (1), \cJ (1) \Bigr\rangle_{\cG(\bv)}}.
\end{equation}
Equation~\eqref{stepsize} mirrors the choice of \citet{Han2014} to use the inner product with the Jacobi field (projecting both sides of the equation onto one side). The inner product used is the one that induces the metric~\eqref{MVNmetric}, evaluated at the geodesic endpoint. Just like the metric, this inner product depends also on the location on the manifold $\btheta$, and for two vectors $\bu = (\bu_{\bmu}, \bu_{\bSigma})$ and $\bu' = (\bu'_{\bmu}, \bu'_{\bSigma})$ takes the form
\begin{equation}
    \label{MVNinnerproduct}
\Bigl\langle \bu, \bu' \Bigr\rangle_{\bp} = \bu_{\bmu}^T \btheta_{\Sigma}^{-1} \bu'_{\bmu} + \frac{1}{2} \tr \Bigl( \btheta_{\bSigma}^{-1} \bu_{\bSigma}  \btheta_{\bSigma}^{-1} \bu'_{\bSigma} \Bigr).
\end{equation}

\begin{figure}
    \centering
    \includegraphics[width=0.85\textwidth, trim={3.5cm, 2cm, 4cm, 3cm},clip]{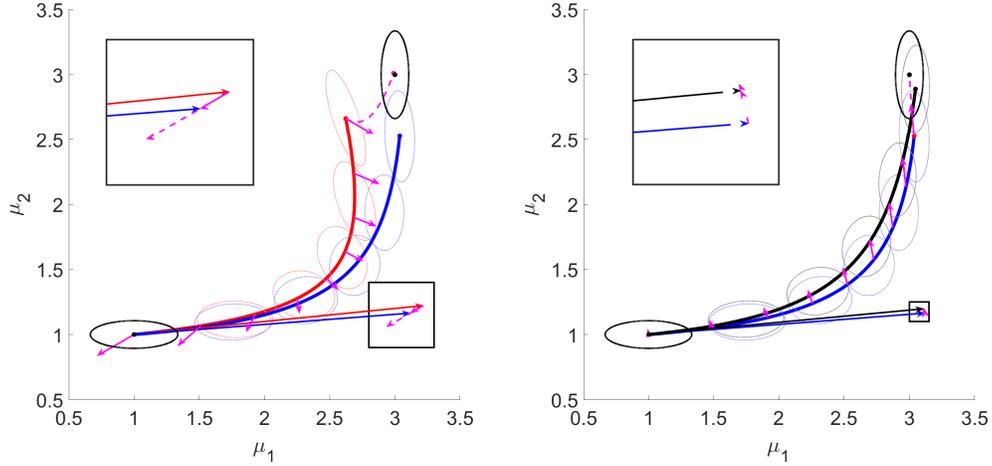}
    \caption{The process of finding geodesics via shooting and parallel transport, visualised for bivariate MVNs. {\bf Left:} A geodesic is fired using an initial velocity guess (red curve). A residual vector (magenta arrow) joining the geodesic's endpoint to the target (black) is found, and parallel transported back along the curve to become a correction to its initial velocity. The update is scaled using equation~\eqref{stepsize} to get the new velocity vector (blue arrow) defining the improved geodesic (blue curve). {\bf Right:} A second shooting iteration improves the geodesic to the black one. Only the $\bmu$-component of each vector is visualised, due to the nature of bivariate MVN visualisation.}
    \label{fig:geodesic_shooting}
\end{figure}

With all of these pieces in place, shooting can proceed by repeatedly iterating the process of choosing a vector connecting the current geodesic endpoint to the actual target, parallel transporting this vector to the origin via backwards integration of equations~\eqref{parallel_transport}, and then choosing a stepsize for this velocity update via equation~\eqref{stepsize} to obtain a new geodesic. Termination can be triggered using the norm of residual vectors, or the symmeterised Kullback--Leibler divergence as an easily-calculated proxy for squared manifold distance,
\[
\DsKL = \frac{1}{4} \Bigl( \tr\bigl(\bSigma_1^{-1}\bSigma_0 + \bSigma_0^{-1} \bSigma_1 - 2\bI\bigr) + (\bmu_1 - \bmu_0)^T \bigl( \bSigma_0^{-1} + \bSigma_1^{-1} \bigr) (\bmu_1 - \bmu_0 ) \Bigr)
\]
Two iterations of the algorithm are demonstrated graphically in~\figref{fig:geodesic_shooting}, and algorithm pseudocode is given in Algorithm~\ref{suppalg:MVNshooting} in the supplement.

\subsection{Alternate Choices of Residual}
\label{sec:improved_shooting}

\begin{figure}
    \centering
    % FIRST FIGURE
    \begin{minipage}{0.015\textwidth}
    {\bf a)}
    \end{minipage}%
    \begin{minipage}{0.475\textwidth}
    \centering
    \includegraphics[width=0.95\textwidth, trim={0cm, 0cm, 2cm, 0cm}, clip]{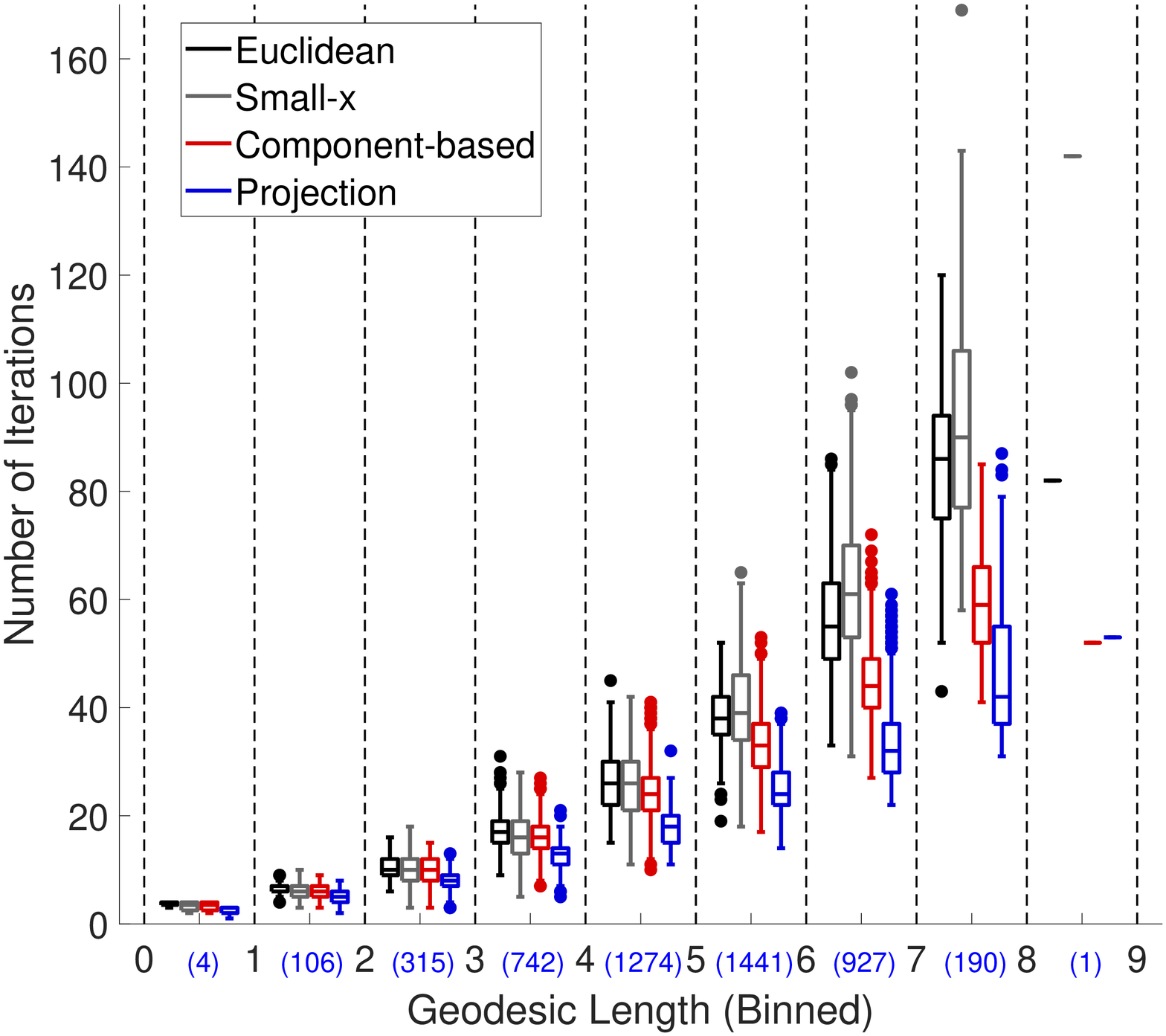}
    \end{minipage}%
    % SECOND FIGURE
    \begin{minipage}{0.015\textwidth}
    {\bf b)}
    \end{minipage}%
    \begin{minipage}{0.475\textwidth}
    \centering
    \includegraphics[width=0.95\textwidth, trim={1cm, 0cm, 2cm, 0cm}, clip]{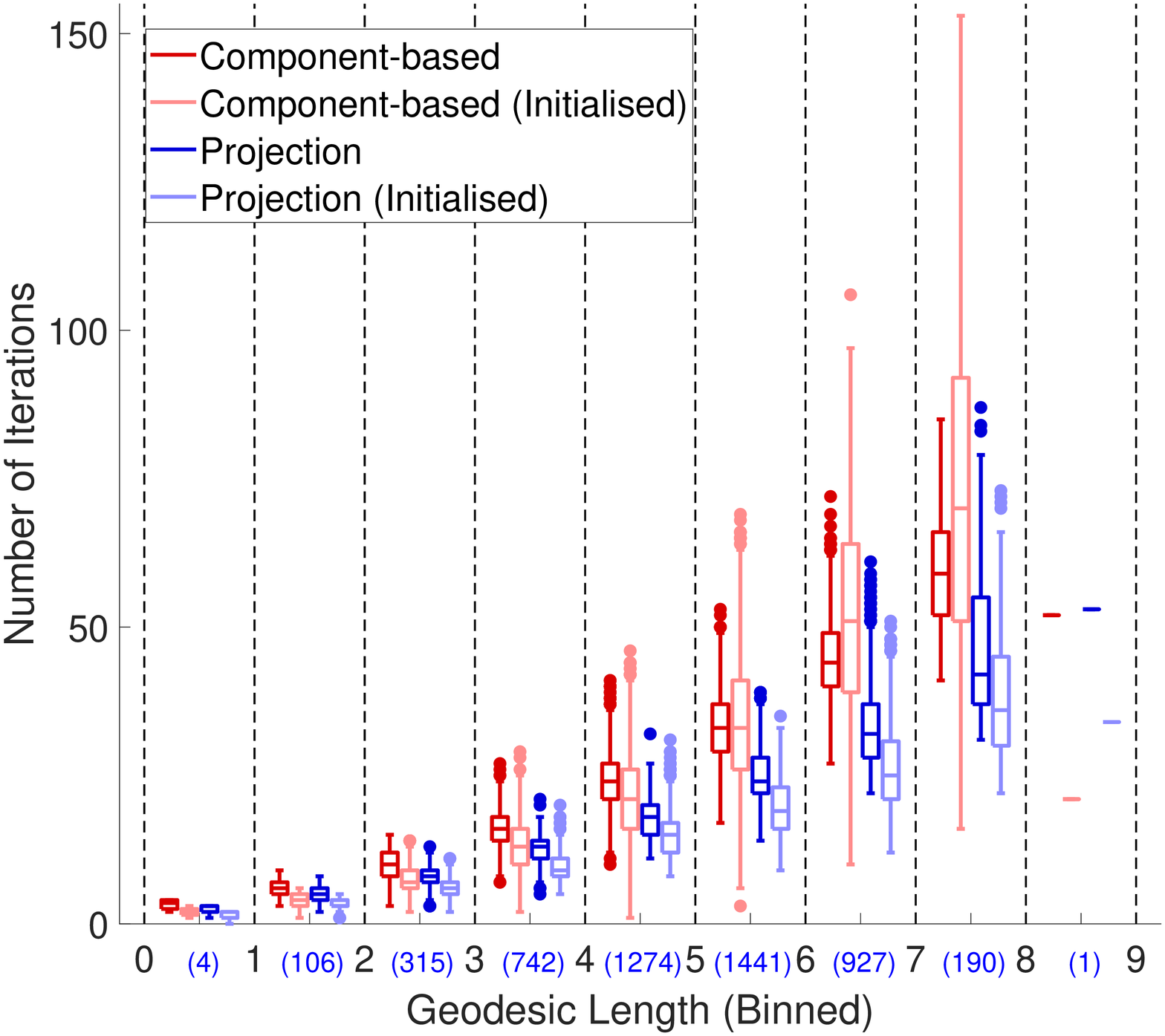}
    \end{minipage}
    \caption{Boxplots demonstrating the convergence rates of shooting for different choices of residual vector over 5000 randomly-generated scenarios. Lengths are binned to aid in summarising the results (numbers in parentheses denote bin sample counts). {\bf a)} Good approximations improve convergence, especially as problem difficulty (distance to target) increases. The projection-based approximation~\eqref{projection_approximation} performs best. {\bf b)} Initialising shooting using a geodesic approximation is generally beneficial, especially using approximation~\eqref{projection_approximation}.}
    \label{fig:residual_choice}
\end{figure}

We find that approximating the geodesic velocity between a pair of points produces significantly better residual vectors for shooting. In~\figref{fig:residual_choice}, we demonstrate iteration counts for shooting across 5000 test problems generated in the same fashion as in Section~\ref{sec:common_paths}. As might be expected, small-$\bx$ approximations begin to hurt convergence as the distance to the target increases and their performance during the initial shooting iterations deteriorates. Defining residuals using the component-based or projection-based approximations consistently and significantly improves mean performance, the latter performing best. 

\citet{Han2014} initialise geodesic shooting using a zero velocity, a choice that we find can help with convergence by avoiding pinning to a location that generates bad shooting updates. However using the geodesic approximations we have reviewed here, we may initialise shooting using a velocity that is already close to the one being sought, potentially improving convergence and reducing the chance of pinning. In~\figref{fig:residual_choice}b we compare between zero-initialised and approximation-initialised shooting, for the two most competitive approximations (component-wise,~\eqref{MVNgeodesics_approx}, and projection-based,~\eqref{projection_approximation}). Pinning was not observed in these test cases, and initialisation generally speeds up performance. However, as distance to the target grows larger and approximation quality decreases, initialisation can become hurtful (observed here for component-based initialisation). Initialisation using the projection-based approximation remains beneficial for the target distances considered here ($d_F < 8$).

\subsection{Shooting via Path Refinement}
\label{sec:path_refinement}

Addressing the issues with convergence of shooting for pairs of MVNs separated by larger Fisher--Rao distances, \citet{Han2014} also proposed a method in which a set of intermediate points between the start point and target are refined until all lie along the desired geodesic path. This is achieved by using the shooting described in Section~\ref{sec:shooting} to find geodesic connections between all of the even-numbered points and all of the odd-numbered points, and using the midpoints of these geodesics to update the locations of the opposite set of points. Iterating this process causes the path defined by these points to settle onto the geodesic (see \figref{fig:path_refine_schematic}, and Algorithm~\ref{suppalg:MVNmultishooting} in the supplementary material), and each individual set (odd or even) may be run in parallel. Once the locations of the set of points have settled onto the geodesic, the stored initial velocity for the geodesic between the first and third points can be scaled according to the number of points to obtain the desired velocity. The Fisher--Rao distance is given simply by totalling the lengths of each individual geodesic segment.

\begin{figure}
    \centering
    \includegraphics[width=0.45\textwidth,trim={1cm, 0cm, 2cm, 1.5cm},clip]{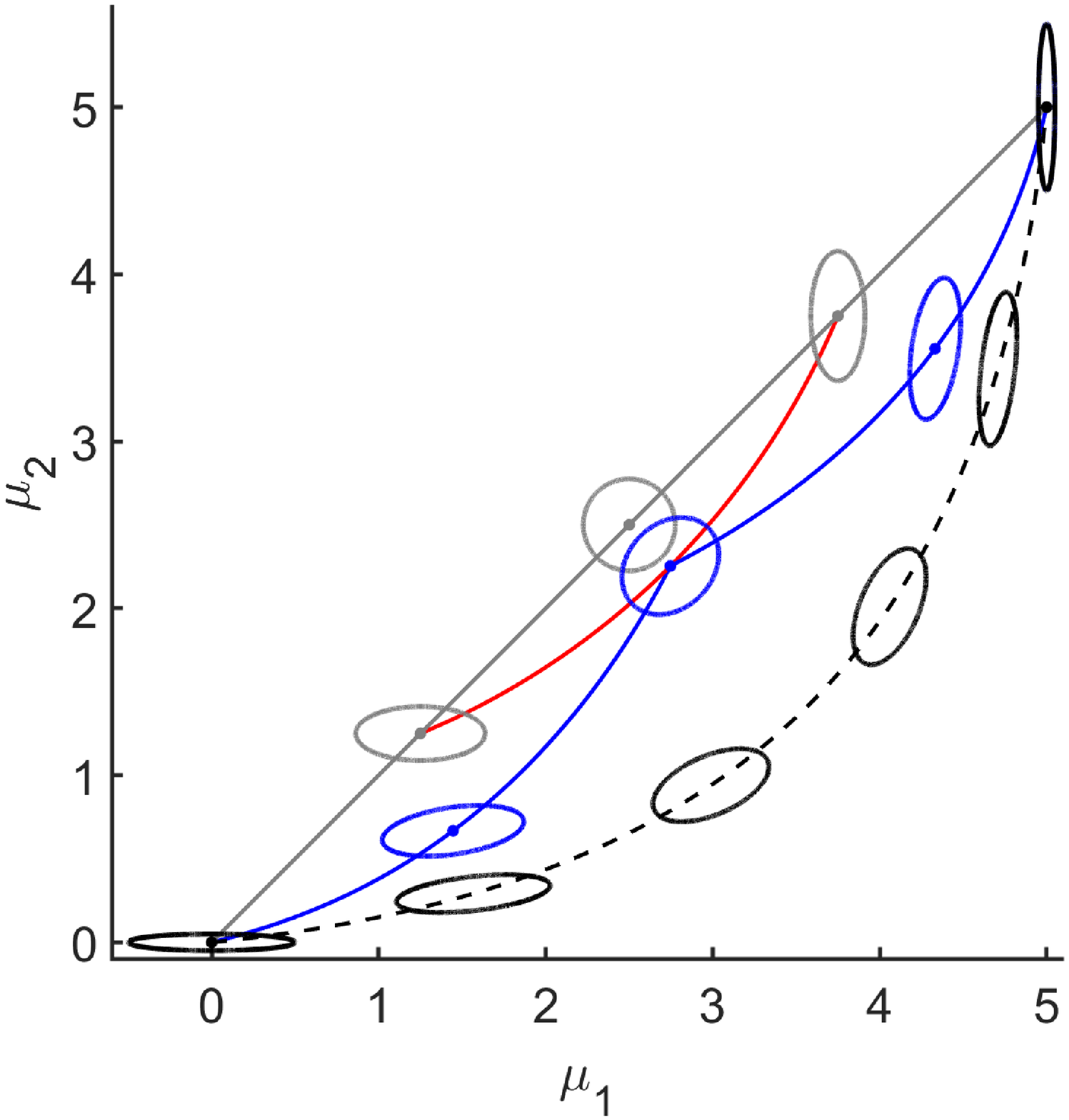} \hspace{0.05\textwidth}
    \includegraphics[width=0.45\textwidth,trim={1cm, 0cm, 2cm, 1.5cm},clip]{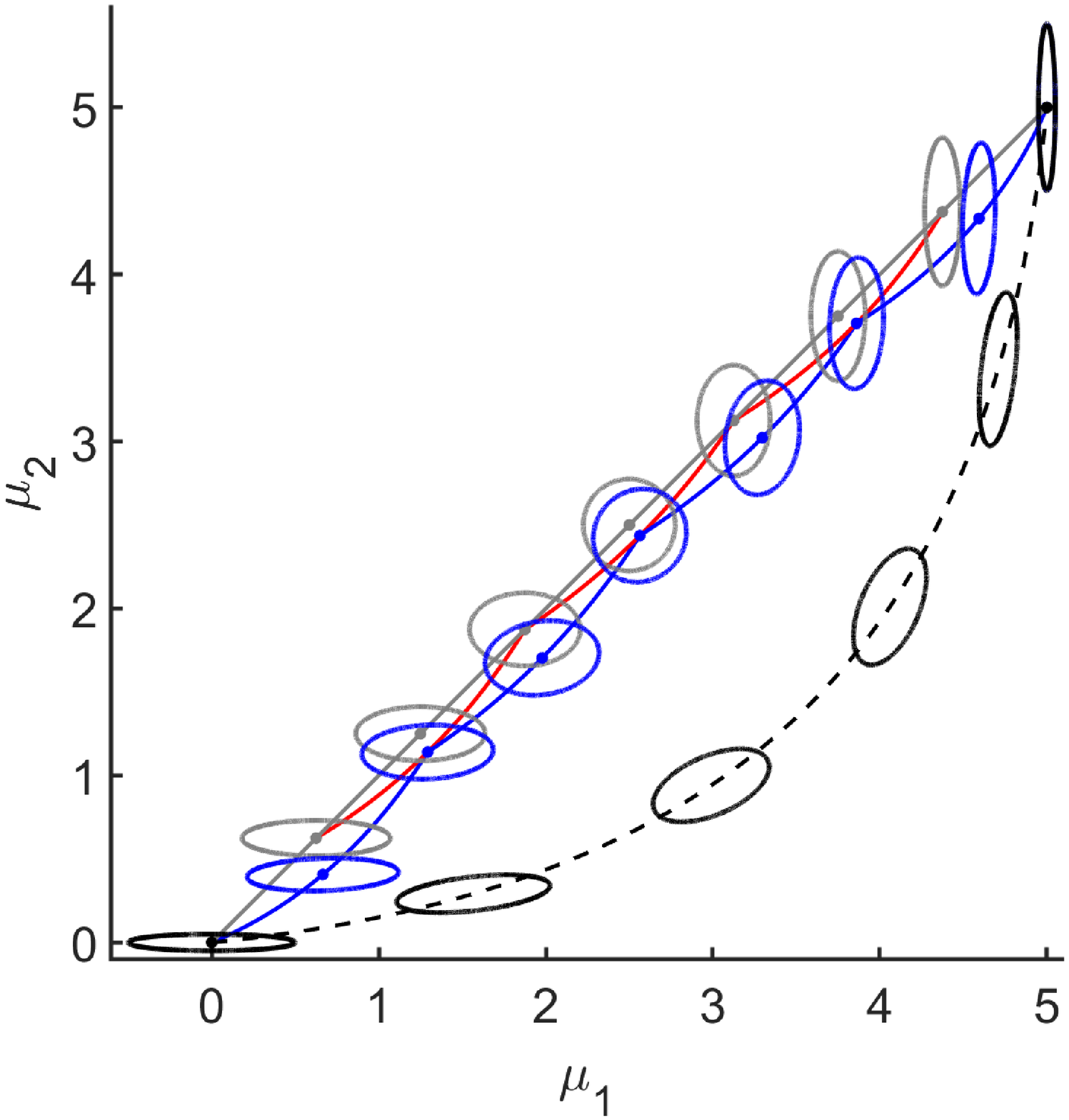}\\[-0.35cm]
    \caption{An example iteration of path refinement, beginning with the Euclidean path (grey) and attempting to find the geodesic path (black, dashed). The geodesic between the even-index points (red) is used to position the central odd-index points. These geodesics (blue), together with their midpoints, provide an updated path that is closer to geodesic. Using more path points (right panel) produces a smaller improvement in path quality per iteration.}
    \label{fig:path_refine_schematic}
\end{figure}

Path refinement breaks down a longer-distance shooting problem into many shorter-distance shooting problems, avoiding the convergence issues with the former. The algorithm requires a full application of the base shooting algorithm (Algorithm~\ref{suppalg:MVNshooting}) per point, {\it per iteration}, although we find this cost can be significantly ameliorated by initialising each shooting call using the geodesic found in the previous path refinement iteration. Instead, the issue with path refinement is that as the number of points defining the path increases, it becomes too easy for individual segments of the path to be approximately geodesic, causing little update in their midpoints (\figref{fig:path_refine_schematic}) while the overall path remains non-geodesic. This also requires us to specify a stricter tolerance for the shooting subproblems than we desire for the overall geodesic being sought.

These effects are demonstrated by using path refinement with a large number of points ($N = 129$) to find the geodesic connecting two distant MVNs, $\left(\begin{bmatrix}1 \\ 2\end{bmatrix}, \begin{bmatrix} 1 & 0.1 \\ 0.1 & 10 \end{bmatrix}\right)$ and $\left(\begin{bmatrix}70 \\ 35\end{bmatrix}, \begin{bmatrix} 10 & -0.8 \\ -0.8 & 1 \end{bmatrix}\right)$ with $d_F = 10.417$. Given their observed superiority in Section~\ref{sec:improved_shooting}, in the results that follow all calls to Algorithm~\ref{suppalg:MVNshooting} as part of the path refinement process use projection-based residuals to formulate shooting updates. However, we find this choice is far less important inside path refinement.

\begin{figure}
    \centering
    \begin{minipage}{0.025\textwidth}%
    \bf{a)} \vspace{2.5cm}
    \end{minipage}%
    \begin{minipage}{0.47\textwidth}%
    \includegraphics[width=1\textwidth, trim={0cm, 0cm, 0cm, 0cm}, clip]{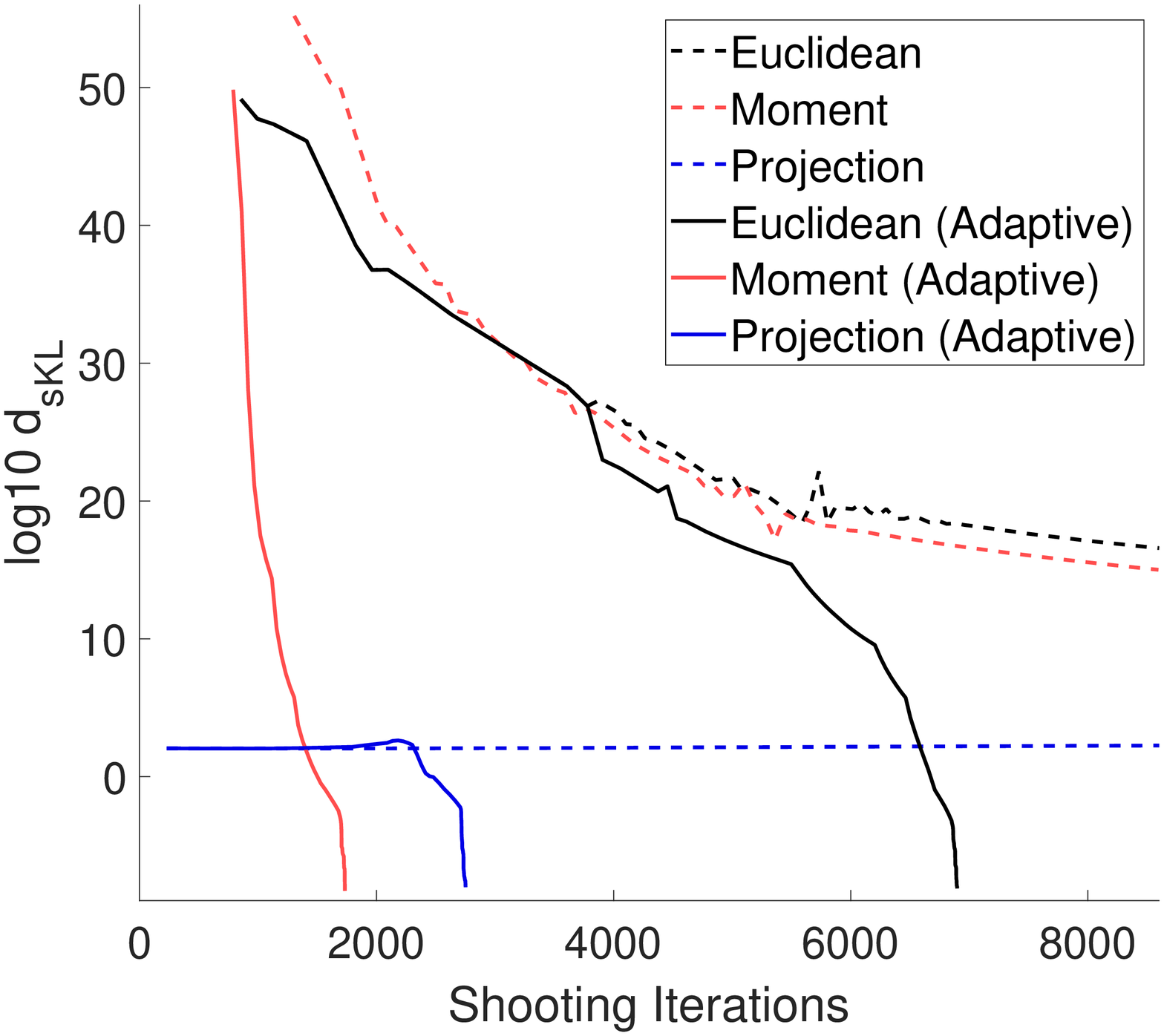}
    \end{minipage}%
    \begin{minipage}{0.035\textwidth}%
    \bf{b)} \vspace{2.5cm}
    \end{minipage}%
    \begin{minipage}{0.47\textwidth}%
    \includegraphics[width=1\textwidth, trim={0cm, 0cm, 0cm, 0cm}, clip]{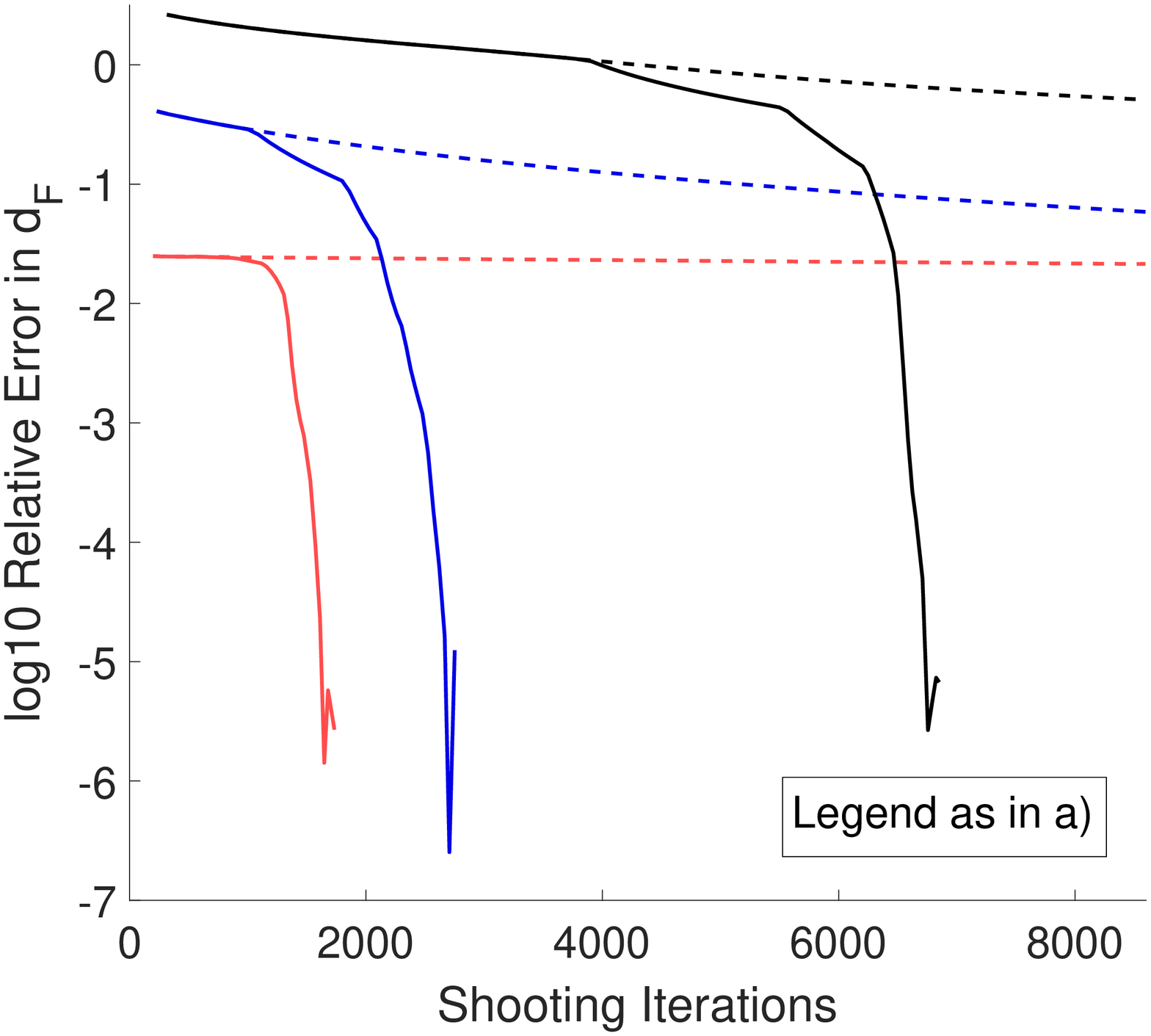}
    \end{minipage}%
    \caption{Performance of standard (dashed lines) and adaptive (solid lines) path refinement for different initialisations (colours), for a large number of path points (N=129) and bivariate MVNs separated by $d_F = 10.42$. {\bf a)} Geodesic velocity accuracy, in terms of distance $\DsKL$ between geodesic endpoint and target. Adaptive path refinement is required for convergence, and initialising with a path closer to geodesic speeds convergence. {\bf b)} Similar performance occurs for Fisher--Rao distance calculation. The standard algorithm can produce false convergence to an incorrect value.}
    \label{fig:path_refine_performance}
\end{figure}

\figref{fig:path_refine_performance}a demonstrates the failure of standard path refinement to converge (dashed lines). After more than 8000 shooting iterations the endpoint of the geodesic using the $\bv$ returned by standard path refinement (with a Euclidean initialisation) is still $D_{sKL} > 10^{15}$ away from the target MVN, highlighting the difficulty of estimating geodesic velocities over large distances (for which small changes have great influence on the endpoint). Performance in terms of Fisher--Rao distance is slightly better (\figref{fig:path_refine_performance}b), but shows a false convergence where subsequent iterations cease to improve an incorrect length estimate. These effects are most pronounced when the number of path points is high, however the failure to reach the tolerance used in subproblems can be observed even for as few as five path points~(\figref{suppfig:path_refine_Npoints}).

Path refinement may be repaired by adaptively culling points as the algorithm proceeds (Algorithm~\ref{suppalg:MVNmultishooting_adaptive}). Initialising path refinement with $2^m + 1$ points allows us to repeatedly delete all even-indexed points whenever the path ceases to sufficiently improve, using the geodesics found using the higher number of points as initialisations to help the shooting problems now between more distant points to converge. Once the path is reduced to only three points, a final run of basic shooting produces the geodesic, using an approximately-geodesic path as a very strong initial guess. This restores the ability of path refinement to converge to the desired geodesic (solid lines in \figref{fig:path_refine_performance}), and to directly specify the desired tolerance (\figref{suppfig:path_refine_Npoints}).

Once using adaptive path refinement, we observe that the paths discussed in Sections~\ref{sec:common_paths} and~\ref{sec:approx_solutions} significantly improve the algorithm's performance. For two means separated by a large distance, as in the test case here, the moment-averaged path is closest to geodesic and provides the fastest convergence. Projection-based paths~\eqref{projection_path} also offer a significant improvement in convergence over the Euclidean initialisations used by~\citet{Han2014}. Although the best choice of path will be problem-specific, it is comparatively cheap to form paths and calculate their lengths before performing any shooting, informing the likely best choice of path to use in joining a given pair of distant MVNs.

\section{Conclusions}

Although they are among the most mathematically tractable of distributional families, MVNs have continued to defy a complete information geometric treatment. Here we have presented a more pragmatic understanding of how the paths of shortest distance along this manifold behave, and how this may inform efforts to determine such a path between two given MVNs. Finding these paths, or their length (the Fisher--Rao distance), is the primary way in which such information geometry is practically applied~\citep{Sim2012, Han2014, Pilte2016, Du2020}.
Despite the non-trivial Riemannian geometry, we have demonstrated how these geodesics permit a more intuitive, traditional geometric understanding that explains when geodesics do or not travel straight lines in $\bmu$-space (and hence, when a closed-form solution is available). This geometric understanding highlights the importance of repeated eigenvalues of the transformed target covariance matrix $\bSigma_0^{-1/2}\bSigma_1\bSigma_0^{-1/2}$, which correspond to more expansive subspaces of points in $\bmu$-space that can be reached via such geodesics. 

Treating the distance along the statistical manifold associated with MVNs as an information distance also provides a powerful means of interpreting commonly-used paths in sequential or tempered sampling methods. Annealing paths select intelligent, variance-informed paths along which to move the location of a distribution, but fail to incorporate additional variance to reduce the cost of this movement. Moment-averaged paths incorporate additional variance to compensate for the naive path of movement for their location, mimicking a feature of geodesic paths and hence explaining their generally lower information cost. Wasserstein-optimal paths completely ignore the informational context.

Without a closed-form solution available, geodesics on the MVN statistical manifold are found via a numerical shooting process. Here we have thoroughly reviewed this process and examined the concerns associated with its practical usage. In particular, we have highlighted the use of geodesic velocity approximations to accelerate shooting, as well as issues with convergence inherent to the shooting via path refinement proposed by \citet{Han2014} and the use of path point culling to resolve them. Considering path information costs also provides good initialisations for the path refinement process.

This work has concentrated on the information geometry of MVNs, acknowledging their ubiquity in statistics and position as a comparatively tractable family that still demonstrates the additional geometric complexities that come with vector-valued random variables. However, many of the ideas we have summarised and discussed here are extensible to other distributional families and contexts. The intuition of temporarily increasing distribution scale in order to lower the information cost associated with changes in location should be applicable to many other distributions, and indeed (univariate) Cauchy and Student's $t$ distributions possess a Fisher geometry very similar to univariate normal distributions~\citep{Mitchell1988}. The moment-averaged paths that also capture this desirable property are well-defined (though not always in closed form) for any exponential family distribution~\citep{Grosse2013}. The MVN-manifold geodesic connecting the Laplace approximations~\citep{Tierney1984} for a pair of distributions might also be used to inform information-efficient paths between the distributions themselves.

A closed-form solution for the geodesic connecting two members of the MVN family has continued to elude researchers in information geometry despite sophisticated treatments of the topic via differential geometry and Lie theory~\citep{Calvo1991, Imai2011}. By presenting here the information geometry of MVNs with only minimal use of differential geometry, our hope is that these fascinating concepts might become more broadly considered, and that this might increase the opportunity for these remaining questions in the space of MVNs to be solved. Until a closed-form solution for the geodesic (and Fisher--Rao distance) between an arbitrary pair of MVNs is found, the numerical approaches we have reviewed and extended here remain necessary.

\section*{Acknowledgements}

BAJL acknowledges infrastructural and research support from the QUT Centre for Data Science. KB is supported by the ARC Centre of Excellence for Plant Success in Nature and Agriculture (CE200100015). RWS is supported by the Coordena\c{c}\~{a}o de Aperfei\c{c}oamento de Pessoal de N\'{i}vel Superior (CAPES), Funda\c{c}\~ao de Amparo \`{a} Pesquisa do Estado de Minas Gerais (FAPEMIG) and Conselho Nacional de Desenvolvimento Cient\'{i}fico e Tecnol\'{o}gico (CNPq).

\section*{Code Availability}

MATLAB implementations of routines for generating, finding and visualising geodesics or other common paths through the family of multivariate normals, including code used to generate the figures in this manuscript, is available via\\ \texttt{https://www.github.com/betalawson/MVNgeodesics}.

\bibliographystyle{unsrtnat}
\bibliography{MVNgeodesics}

\clearpage

% Initialise appendix
\appendix

% Relabel appendix sections with S
\renewcommand{\thesection}{S}

% Reset figure counter, and label supplemental figures as S1, S2, ...
\setcounter{figure}{0}
\makeatletter 
\renewcommand{\thefigure}{S\@arabic\c@figure}
\makeatother

% Relabel algorithms as S1, S2, ...
\makeatletter 
\renewcommand{\thealgorithm}{S\@arabic\c@algorithm}
\makeatother

\section*{Supplementary Material}

\subsection{Derivation of the ``Axis-aligned'' Solution}

For a geodesic beginning at the origin, the geodesic equations for multivariate normals take the form~\eqref{MVNgeodesic_solution}. Using the affine transformation to shift the geodesic between two points to one beginning at the origin, the target point is given by equation~\eqref{affine_origin}, or in terms of canonical co-ordinates for MVNs,
\[
\bdelta_t = \bR^T \bSigma_0^{-1/2} \bSigma_1^{-1} \left(\bmu_1 - \bmu_0 \right), \qquad \qquad \bDelta_t = \bR^T \bSigma_0^{-1/2} \bSigma_1 \bSigma_0^{-1/2} \bR.
\]
The remaining problem for finding the geodesic is to identify the initial velocity parameters $\bx$ and $\bB$, such that at $t = 1$ the geodesic passes through the target point,
\[
\begin{aligned}
        \bdelta_t &= - \bB \Bigl( \cosh(\bG) - \bI \Bigr) (\bG^{-})^{2}\bx + \sinh(\bG ) \bG^{-} \bx, \\
        \bDelta_t &= \bI + \frac{1}{2} \Bigl( \cosh(\bG ) - \bI  \Bigr) + \frac{1}{2} \bB \Bigl( \cosh(\bG) - \bI \Bigr) (\bG^{-})^{2} \bB - \frac{1}{2} \sinh(\bG) \bG^{-} \bB - \frac{1}{2} \bB \sinh(\bG) \bG^{-}.
\end{aligned}
\]

As discussed in the main text, these solutions occur when $\bB$ is a diagonal matrix, and $\bx$ is aligned with a single co-ordinate axis (which also ensures $\bG$ is diagonal). Considering the behaviour of geodesics in equation~\eqref{MVNgeodesic_solution} under these conditions, we see that $\bdelta_t$ remains aligned with the same axis as $\bx$, and $\bDelta$ remains a diagonal matrix. If the target point falls within this space of valid targets, then we have a solution to the geodesic connecting the origin to that target. In this case where $\bB$ and $\bG$ are diagonal, the equations decouple in the sense that they may be solved component-wise,
\[
\begin{aligned}
        \delta_k &= \left( - \frac{B_{kk}}{G_{kk}^2} \left( \cosh G_{kk} - 1 \right)  + \frac{1}{G_{kk}} \sinh G_{kk} \right) x_k, \\
        \Delta_{kk} &= 1 + \frac{1}{2} \left( \cosh G_{kk} - 1 \right) + \frac{1}{2} \frac{B_{kk}^2}{G_{kk}^2} ( \cosh G_{kk} - 1 ) - \frac{B_{kk}}{G_{kk}} \sinh G_{kk}.
\end{aligned}
\]

Each pair of equations can be related to the geodesics of the Poincar{\'{e}} half-plane model~\citep{Skovgaard1984}. For completeness, we present here the direct approach given in~\citep{Sim2012} in a little extra detail. For notational simplicity, we drop the component notation in what follows, such that non-bolded quantities refer to the value in one component (diagonal component for matrices).

Defining $R = B / G$, the second equation can be rewritten
\[
\Delta = \frac{1}{2}\left( 1 - R^2 \right) + \frac{1}{2} \left( 1 + R^2 \right) \cosh G - R \sinh G.
\]
Substitution of the exponential definitions for $\cosh$ and $\sinh$ gives
\[
\begin{aligned}
\Delta &= \frac{1}{2}\left( 1 - R^2 \right) + \frac{1}{4} \left( 1 + R^2 \right) \left(  e^{G} + e^{-G} \right) - \frac{1}{2} R \left( e^{G} - e^{-G} \right) \\
            &= \frac{1}{2} \left(1 - R^2 \right) + \left( \frac{1}{4} + \frac{1}{4} R^2 - \frac{1}{2} R \right) e^{G} + \left( \frac{1}{4} + \frac{1}{4} R^2 + \frac{1}{2} R \right) e^{-G} \\
            &= \frac{1}{2} \left(1 - R^2 \right) + \frac{1}{4} \Bigl( 1 - R \Bigr)^2  e^{G} + \frac{1}{4} \Bigl( 1 + R \Bigr)^2  e^{-G} \\
            &= \frac{1}{4} \Bigl( (1 - R) e^{\frac{1}{2} G} + (1 + R) e^{-\frac{1}{2} G} \Bigr)^2.
\end{aligned}
\]
Given that $R^2 = \dfrac{B^2}{G^2} = \dfrac{B^2}{B^2 + 2x^2} \leq 1$, the right hand side equation is a positive term squared, and as such we can take the positive branch when taking the square root of both sides. This then gives a quadratic in $e^{\frac{1}{2} G}$, with solution
\[
\begin{aligned}
    e^{\frac{1}{2}G} &= \frac{\sqrt{\Delta} \pm \sqrt{ \Delta - 1 + R^2}}{(1-R)}.
\end{aligned}
\]
Given that $R \leq 1$, both branches give positive values and are potentially valid. Re-arranging for $G$ gives
\[
G = \log \left( \frac{\sqrt{\Delta} \pm \sqrt{ \Delta - 1 + R^2}}{1-R} \right)^2.
\]

Now returning to the first equation of the component-wise solution and making the same substitution for $R$ including $x = \pm G \sqrt{\frac{1}{2}(1-R^2)}$, we have
\[
\begin{aligned}
\delta &= \Bigl( - \frac{B}{G^2} ( \cosh G - 1 )  + \frac{1}{G} \sinh G \Bigr) x \\
         &= \pm \Bigl( - R \cosh G + R + \sinh G \Bigr)  \sqrt{\frac{1}{2}(1 - R^2)}
\end{aligned}
\]
Setting $\nu = \left( \frac{\sqrt{\Delta} \pm \sqrt{ \Delta - 1 + R^2}}{(1-R)} \right)^2$ such that $G = \log(\nu)$ for notational convenience, we have that
\[
\begin{aligned}
\delta &= \pm \Biggl( - R \left( \frac{\nu}{2} + \frac{1}{2\nu} \right) + R + \frac{\nu}{2} - \frac{1}{2\nu} \Biggr) \sqrt{\frac{1}{2}(1 - R^2)} \\
\delta^2 &= \frac{1}{2} \Biggl( - R \left( \frac{\nu}{2} + \frac{1}{2\nu} \right) + R  + \frac{\nu}{2} - \frac{1}{2\nu} \Biggr)^2 \left(1 - R^2 \right) \\
          &= \frac{1}{8\nu^2} \left( -R \nu^2 - R + 2 R \nu + \nu^2 - 1 \right)^2 \left(1 - R^2 \right) \\
          &=\frac{1}{8\nu^2}(\nu-1)^2 \Bigl( \nu (1-R) + R + 1 \Bigr)^2 (1-R^2).
\end{aligned}
\]
Substituting the definition of $\nu$, and after some simplification,
\[
\begin{aligned}
\delta^2 &= \frac{\Bigl(\Delta \pm \sqrt{\Delta}\sqrt{\Delta + 1 - R^2} - (1-R)\Bigr)^2\Bigl(\Delta \pm \sqrt{\Delta}\sqrt{\Delta + 1 - R^2}\Bigr)^2}{2 \Bigl(\Delta \pm \sqrt{\Delta} \sqrt{\Delta - 1 + R^2} + \frac{1}{2}\bigl(1 - R^2\bigr)\Bigr)^2} \frac{1+R}{1-R}.
\end{aligned}
\]
From here, computer algebra gives the solution
\[
R = \frac{\pm(\delta^2 + 2\Delta - 2\Delta^2)}{\sqrt{4\Delta^4 + 4\Delta^2\delta^2 + \delta^4 - 8\Delta^3 + 4\Delta\delta^2 + 4\Delta^2}}.
\]
Note that this is true regardless of whether the positive or negative branch in the definition of $\nu$ is selected.

As $G$ is defined only in terms of $G^2$ (consider the Taylor series representation of the equations for $\delta$ and $\Delta$), we select $G$ to always be positive. Then, considering that in the scenario that the mean does not need to move ($\delta = 0$), a decrease in covariance ($\Delta > 1$) must correspond to a negative value for $B$, we know to choose the positive branch. As such, using the notation-simplifying variables from \tabref{tab:known_geodesics} in the main document, we write
\[
R = \frac{\alpha}{\gamma},
\]
and comparing the expansion of $\alpha^2$ to the terms in $\gamma^2$ reveals the relation
\[
\gamma = \sqrt{ \alpha^2 + 8\Delta^2 \delta^2}.
\]

Substituting the definition of $R$ into the definition for $\nu$, and using the relation $G = \log(\nu) = \acosh\left(1 + \dfrac{(\nu-1)^2}{2 \nu}\right)$, it can be shown that for one of the two branches of $\nu$,
\[
G = \acosh\left(1 + \frac{\gamma^2}{8\Delta^3}\right),
\]
which matches the form for $g$ given in~\tabref{tab:known_geodesics}. Then, with $G$ found, we have that
\[
B = R G = \frac{\alpha}{\gamma} G,
\]
and for $x$,
\[
\begin{aligned}
x &= \pm G \sqrt{ \frac{1}{2} \Biggl( 1 - \left(\frac{\alpha}{\gamma}\right)^2 \Biggr) }\\
  &= \pm G \sqrt{ \frac{1}{2} \cdot 8 \frac{\Delta^2 \delta^2}{\gamma^2} } \\
  &= \frac{2 G \delta \Delta}{\gamma}.
\end{aligned}
\]
The positive branch is selected so that the sign of $x$ matches the sign of $\delta$, such that $\mu$ moves in the correct direction.

In the above, component notation was dropped for simplicity, and the above procedure must be solved for each component separately. However, this solution is only valid when $\bB$ and $\bG$ are diagonal matrices, corresponding to the case where $\bdelta = \delta \be_k$ is aligned with a single co-ordinate axis. As such, only one component takes the full form elucidated above. In the remainder of cases, $\delta = 0$, which gives $x = 0$ and thus $B = G$ and $R = 1$. The second equation for the form of the geodesics then reduces simply to
\[
\begin{aligned}
\Delta &= \cosh B - \sinh B\\
       &= e^{-B} \\
\log(\Delta^{-1}) &= B \\
\log(\Sigma_t) &= B,
\end{aligned}
\]
where $\Sigma_t$ is the target covariance for the current component. Using this velocity in the geodesic equation gives
\[
\Sigma(t) = e^{t \log \Sigma_t} = \left(\Sigma_t\right)^t,
\]
showing that the covariance transforms exponentially to its target value, starting from the origin $\Sigma = 1$ as we are in the transformed space. This corresponds to the observation in the main text that the covariance eigenvalues transform exponentially in directions other than the direction of movement --- that result is obtained simply by removing the diagonalising part of the affine transform.

Under these conditions where $\bB$ and $\bG$ are diagonal and the geodesic equations can be considered component wise, we obtain our solutions for $\bx$ and $\bB$ by now adding these components,
\[
\bx = \frac{2 G \delta \Delta}{\gamma} \be_k, \qquad \qquad \bB = \frac{\alpha}{\gamma} G \be_{k} \be_{k}^T - \sum_{j \neq k} \log D_{jj} \be_{j} \be_{j}^T,
\]
where the index $k$ is defined by the condition $\bdelta = \delta \be_{k}$, that is, the component in which the mean must move to reach the target.

\subsection{A Small-$\bx$ Velocity Approximation via Taylor Series}

When $\bx$ is small, $\bG^2 = \bB^2 + 2\bx \bx^T \approx \bB^2$ and equations~\eqref{MVNgeodesic_solution} for the form of the geodesics simplify considerably. The second equation becomes
\[
\bDelta(t) \approx \bI + \frac{1}{2} \Bigl(\cosh(t\bB) - \bI \Bigr) + \frac{1}{2} \Bigl( \cosh(t\bB) - \bI \Bigr)  - \frac{1}{2}\sinh(t\bB) - \frac{1}{2}\sinh(t\bB) = e^{-t\bB},
\]
and treating this as an equality and choosing $\bB$ to satisfy $\bDelta(1) = \bDelta_t$ gives the approximation $\bB \approx -\log(\bDelta_t) = \log(\bSigma_t)$. Expanding the Taylor series for $\cosh$ and $\sinh$ in the first of equations~\eqref{MVNgeodesic_solution}, one obtains
\[
\bdelta(t) = -\bB \sum_{j=1}^{\infty} \frac{1}{2j!} (t \bG)^{2j-2} \bx + \sum_{j=1}^{\infty} \frac{1}{(2j-1)!} (t \bG)^{2j-2} \bx,
\]
which after substituting in $t=1$ and again applying $\bG^2 \approx \bB^2$ gives
\[
\begin{aligned}
\delta_t &\approx \Bigl( - \sum_{j=1}^{\infty} \frac{1}{2j!} \bB^{2j-1} + \sum_{j=1}^{\infty} \frac{1}{(2j-1)!} \bB^{2j-2} \Bigr) \bx \\
          &\approx \bB^{-1} \Bigl( - \sum_{j=1}^{\infty} (-1)^j \frac{1}{j!} \bB^{j} \Bigr) \bx \\
          &\approx \bB^{-1} \Bigl( \bI - e^{-\bB} \Bigr) \bx.
\end{aligned}
\]
Substituting the approximation for $\bB$, $\bB \approx - \log \bDelta_t$, this gives for $\bx$
\[
\bx \approx \Bigl( \bI - \bDelta \Bigr)^{-1}  - \log (\bDelta_t) \bdelta_t.
\]
Expressing these approximate choices for $\bx$ and $\bB$ in terms of the mean and covariance parameters of the target multivariate normal distribution, they become
\[
\begin{aligned}
\bx &\approx \Bigl( \bI - (\bI - \bSigma_t)^{-1} \Bigr) \log( \bSigma_t) \bSigma^{-1}_t \bmu_t\\
\bB &\approx \log \bSigma_t,
\end{aligned}
\]
where the Sherman--Morrison--Woodbury formula is used to write
\[
\Bigl( \bI - \bA^{-1} \Bigr)^{-1} = \bI - (\bI - \bA)^{-1}.
\]

\subsection{A Velocity Approximation via Higher-dimensional Manifold Embedding}

The approximately geodesic paths found by embedding the MVN statistical manifold in the space of symmetric positive definite matrices also produce a velocity approximation, simply by taking the initial velocity of these paths. As discussed in the main document, these paths take the form
\[
\bp(t) = \begin{bmatrix} \bSigma_t + \bmu_t \bmu_t^T & \bmu_t \\ \bmu_t^T & 1 \end{bmatrix}^{t}.
\]
We first seek $\bp'(0)$, which is found by taking the logarithm of both sides and differentiating,
\[
\deb{t} \log \bp(t) = \log \begin{bmatrix} \bSigma_t + \bmu_t \bmu_t^T & \bmu_t \\ \bmu_t^T & 1 \end{bmatrix}
\]
Then, as for the origin-transformed problem we have $\bp(0) = \bI$, $\bp(0)$ and $\bp'(0)$ commute. As such we may use $\deb{t} \left( \log \bp \right) = \bp^{-1} \de{\bp}{t}$ to obtain simply
\[
\bp'(0) = \log \begin{bmatrix} \bSigma_t + \bmu_t \bmu_t^T & \bmu_t \\ \bmu_t^T & 1 \end{bmatrix}.
\]

Now we take the derivative of the embedding of an MVN path $\btheta(t) = \bigl(\bmu(t), \bSigma(t)\bigr)$ into the same space,
\[
\begin{aligned}
\bp_{\cN}'(t) &= \deb{t} \begin{bmatrix} \bSigma(t) + \bmu(t) \bigl[\bmu(t)\bigr]^T & \bmu(t) \\ \bigl[\bmu(t)\bigr]^T & 1 \end{bmatrix} \\
&= \begin{bmatrix} \bSigma'(t) + \bmu'(t) \bigl[\bmu(t)\bigr]^T + \bmu(t) \bigl[\bmu'(t)\bigr]^T & \bmu'(t) \\ \bigl[\bmu'(t)\bigr]^T & 0 \end{bmatrix}.
\end{aligned}
\]
At $t = 0$, $\bmu(0) = \b0$, $\bSigma(0) = \bI$ and the initial velocity defines $\bmu'(0) = \bx$ and $\bSigma'(0) = \bB$,
\[
\bp'_{\cN}(0) = \begin{bmatrix} \bB & \bx \\ \bx^T & 0 \end{bmatrix}.
\]
Equating the two then provides the approximate choice of $\bx$ and $\bB$,
\[
\begin{bmatrix} \bB & \bx \\ \bx^T & 0 \end{bmatrix} \approx \log \begin{bmatrix} \bSigma_t + \bmu_t \bmu_t^T & \bmu_t \\ \bmu_t^T & 1 \end{bmatrix},
\]
where $\bx$ and $\bB$ can be read out from the appropriate matrix blocks on the right-hand side. When $\bmu_t = \bx = \b0$, this approximation reduces to $\bB = \log \bSigma_t$ and hence recovers the ``equal means'' solution exactly (see \tabref{tab:known_geodesics} in main document).

\clearpage

\subsection{Algorithms}
%%% EXTRA ALGORITHM SETUP
\algnewcommand\algorithmicswitch{\textbf{switch}}%
\algnewcommand\algorithmiccase{\textbf{case}}%
\algdef{SE}[SWITCH]{Switch}{EndSwitch}[1]{\algorithmicswitch\ #1\ \algorithmicdo}{\algorithmicend\ \algorithmicswitch}%
\algdef{SE}[CASE]{Case}{EndCase}[1]{\algorithmiccase\ #1}{\algorithmicend\ \algorithmiccase}%
\algtext*{EndSwitch}%
\algtext*{EndCase}%

\begin{algorithm}
\begin{algorithmic}
\Procedure{MVNshooting}{$\btheta, \btheta'$, $\bv_0$, \texttt{tol}, \texttt{r\_norm\_max}}
\State $\bO \gets \left( \b0, \quad \bI \right)$    \Comment{Define origin}
\State $\btheta_t \gets \left(\btheta_{\bSigma}^{-1/2} (\btheta'_{\bmu} - \btheta_{\bmu}),\quad \btheta_{\bSigma}^{-1/2} \btheta'_{\bSigma} \btheta_{\bSigma}^{-1/2}\right)$    \Comment{Target point after $\btheta$ is mapped to origin}
\State $\bv \gets \bv_0$    \Comment{Initialise velocity to provided guess (can use zero)}
\State $\bg \gets \cG(\bv)$    \Comment{Endpoint $\bg$ for the initial geodesic via equation~\eqref{MVNgeodesic_solution_matrix}}
\While{ \Call{symKL}{$\bg, \btheta_t$} $>$ \texttt{tol} }
\State $\br \gets$ \Call{residualVector}{$\bg, \bp_t$}    \Comment{Residual between geodesic endpoint and target}
\State \texttt{r\_norm} $\gets \Bigl[\Bigr.\Call{innerProduct}{\br, \br, \bg}\Bigl.\Bigr]^{1/2}$    \Comment{Used for stepsize considerations}
\State $\Delta \bv \gets$ \Call{backTransport}{$\br, \bv$}    \Comment{Convert residual to velocity correction direction}
\State $h \gets $ \texttt{tol} / $\Bigl[\Call{innerProduct}{\Delta \bv, \Delta \bv, \bO}\Bigr]^{1/2}$ \Comment{Select finite differencing stepsize}
\State $\bj \gets \bigl( \cG(\bv + h \Delta \bv) - \bg \bigr) / h$    \Comment{Evaluate geodesic endpoint sensitivity as per \eqref{jacobi_field}}
\State $s \gets $ \Call{innerProduct}{$\btheta_t - \bg, \bj, \bg$} / \Call{innerProduct}{$\bj, \bj, \bg$}    \Comment{Step length}
\State $s \gets $ \Call{min}{$s$,\texttt{r\_norm\_max} / \texttt{r\_norm} $s$} \Comment{Optionally avoid over-ambitious updates}
\State $\bv \gets \bv + s \Delta \bv$    \Comment{Apply velocity update}
\State $\bg \gets \cG(\bv)$    \Comment{Endpopint of geodesic with updated velocity via~\eqref{MVNgeodesic_solution_matrix}}
\EndWhile
\State $\bv \gets \left(\btheta_{\bSigma}^{1/2} \bvmu, \quad \btheta_{\bSigma}^{1/2} \bvSigma \btheta_{\bSigma}^{1/2}\right)$    \Comment{Transform velocity back into original space}
\State \Return $\bv$     \Comment{And/or return Fisher--Rao distance, $\Bigl[\Call{innerProduct}{\bv, \bv, \btheta}\Bigr]^{1/2}$}
\EndProcedure
\end{algorithmic}
\caption{Procedure for finding geodesics between MVNs $\btheta$ and $\btheta'$ via shooting.}
\label{suppalg:MVNshooting}
\end{algorithm}

\begin{algorithm}
\begin{algorithmic}
\Procedure{MVNpathRefinement}{$\btheta, \btheta'$, \texttt{tol}, \texttt{r\_norm\_max}, N}

%\State $\bp_t \gets (\bSigma_0^{-1/2} (\bmu_1 - \bmu_0), \bSigma_0^{-1/2} \bSigma_1 \bSigma_0^{-1/2})$    \Comment{Target point $\bp_1$ after affine transforming geodesic to origin}
\State $(\ba_1, \ba_2, \ldots, \ba_n) \gets $ \Call{closedFormPath}{$\btheta, \btheta', N$}    \Comment{$N$-point path between $\btheta$ and $\btheta'$}
\State \texttt{looping} $\gets$ \texttt{true}  
\While{\texttt{looping}}
\For{$i \gets 3$ by $2$ to $N-1$}    \Comment{Connect even-index pts. to update odd-index pts.}
    \State $\bv_{i-1} \gets$ \Call{MVNshooting}{$\ba_{i-1}, \ba_{i+1}, \bv_{i-1}$, \texttt{tol}, \texttt{r\_norm\_max}}    \Comment{Find geodesic}
    \State $\ba_{i} \gets \cG_{\ba_{i-1}}(0.5\bv_{i-1})$    \Comment{Place point halfway along geodesic}
    \State $L_{i-1} \gets \Bigl[\Call{innerProduct}{0.5\bv_{i-1},0.5\bv_{i-1},\ba_{i-1}}\Bigr]^{1/2}$    \Comment{Store $d_F(\ba_{i-1},\ba_{i})$}
\EndFor

\For{$i \gets 2$ by $2$ to $N-1$}    \Comment{Connect odd-index pts. to update even-index pts.}
    \State $\bv_{i-1} \gets$ \Call{MVNshooting}{$\ba_{i-1}, \ba_{i+1}, \bv_{i-1}$, \texttt{tol}, \texttt{r\_norm\_max}}    \Comment{Find geodesic}
    \State $\ba_{i} \gets \cG_{\ba_{i-1}}(0.5\bv_{i-1})$    \Comment{Place point halfway along geodesic}
    \State $L_{i-1} \gets \Bigl[\Call{innerProduct}{0.5\bv_{i-1},0.5\bv_{i-1},\ba_{i-1}}\Bigr]^{1/2}$    \Comment{Store $d_F(\ba_{i-1},\ba_{i})$}
\EndFor
\State \texttt{Ltot} $\gets $ \Call{sum}{$L_i$}    \Comment{Calculate current manifold distance of the path}
\State $\bg \gets \cG_{\btheta}(0.5(N-1) \bv_1)$    \Comment{Fire geodesic from first point out to target}
\If{ $\Call{symKL}{\bg, \bq} < $ \texttt{tol}}    \Comment{Terminate loop if target point reached}
    \State \texttt{looping} $\gets$ \texttt{false}
\EndIf
\EndWhile
\State \Return $0.5(N-1) \bv_{1}$ \Comment{And/or return Fisher--Rao distance, \texttt{Ltot}}

\EndProcedure
\\ \\
{\it Note the use here of notation $\cG_{\btheta}(\bv)$ refers to tracing out the geodesic with initial velocity $\bv$ starting from point $\btheta$ (by mapping to the origin, using equation~\eqref{MVNgeodesic_solution_matrix}, and then mapping back). The base notation $\cG(\bv)$ is equivalent to $\cG_{(\b0,\bI)}(\bv)$.}

\end{algorithmic}
\caption{The path refinement process for finding the geodesic between $\btheta$ and $\btheta'$}
\label{suppalg:MVNmultishooting}
\end{algorithm}

\begin{algorithm}
\begin{algorithmic}
\Procedure{MVNadaptivePathRefinement}{$\btheta, \btheta'$, \texttt{tol}, \texttt{r\_norm\_max}, $m$}

%\State $\bp_t \gets (\bSigma_0^{-1/2} (\bmu_1 - \bmu_0), \bSigma_0^{-1/2} \bSigma_1 \bSigma_0^{-1/2})$    \Comment{Target point $\bp_1$ after affine transforming geodesic to origin}
\State \texttt{req\_improvement} $ \gets 0.01$    \Comment{Cull if an iteration fails to reduce relative length by 1\%}
\State $N \gets 2^{m+1}$    \Comment{Number of points for level $l$}
\State $(\ba_1, \ba_2, \ldots, \ba_n) \gets $ \Call{closedFormPath}{$\btheta, \btheta', N$}    \Comment{$N$-point path between $\btheta$ and $\btheta'$}
\State \texttt{looping} $\gets$ \texttt{true}
\State \texttt{Ltot\_old} $\gets \infty$
\While{\texttt{looping}}
\For{$i \gets 3$ by $2$ to $N-1$}    \Comment{Connect even-index pts. to update odd-index pts.}
    \State $\bv_{i-1} \gets$ \Call{MVNshooting}{$\ba_{i-1}, \ba_{i+1}, \bv_{i-1}$, \texttt{tol}, \texttt{r\_norm\_max}}    \Comment{Find geodesic}
    \State $\ba_{i} \gets \cG_{\ba_{i-1}}(0.5\bv_{i-1})$    \Comment{Place point halfway along geodesic}
    \State $L_{i-1} \gets \Bigl[\Call{innerProduct}{0.5\bv_{i-1},0.5\bv_{i-1},\ba_{i-1}}\Bigr]^{1/2}$    \Comment{Store $d_F(\ba_{i-1},\ba_{i})$}
\EndFor

\For{$i \gets 2$ by $2$ to $N-1$}    \Comment{Connect odd-index pts. to update even-index pts.}
    \State $\bv_{i-1} \gets$ \Call{MVNshooting}{$\ba_{i-1}, \ba_{i+1}, \bv_{i-1}$, \texttt{tol}, \texttt{r\_norm\_max}}    \Comment{Find geodesic}
    \State $\ba_{i} \gets \cG_{\ba_{i-1}}(0.5\bv_{i-1})$    \Comment{Place point halfway along geodesic}
    \State $L_{i-1} \gets \Bigl[\Call{innerProduct}{0.5\bv_{i-1},0.5\bv_{i-1},\ba_{i-1}}\Bigr]^{1/2}$    \Comment{Store $d_F(\ba_{i-1},\ba_{i})$}
\EndFor
\State \texttt{Ltot} $\gets $ \Call{sum}{$L_i$}    \Comment{Calculate current manifold distance of the path}
\State $\bg \gets \cG_{\btheta}(0.5(N-1) \bv_1)$    \Comment{Fire geodesic from first point out to target}
\If{ $\Call{symKL}{\bg, \bq} < $ \texttt{tol}}    \Comment{Terminate loop if target point reached}
    \State \texttt{looping} $\gets$ \texttt{false}
\EndIf
\If{\texttt{Ltot} / \texttt{Ltot\_old} $> 1 - \texttt{req\_improvement}$ \& $N > 3$}
\State $N \gets (N-1)/2 + 1$    \Comment{Delete even-index points}
\For{$i \gets 1$ to $N$}    
\State $\ba_i \gets \ba_{2i-1}$    \Comment{Locations of retained points}
\State $\bv_i \gets 2 \bv_{2i-1}$    \Comment{Retained geodesics (doubled as each goes twice as far)}
\EndFor
\EndIf
\EndWhile
\State $\bv \gets \Call{MVNshooting}{\btheta, \btheta', \bv_1, \texttt{tol}, \texttt{r\_norm\_max}}$    \Comment{Perform final shooting}
\State \Return $\bv$ \Comment{And/or return Fisher--Rao distance, $\Bigl[\Call{innerProduct}{\bv,\bv,\btheta}\Bigr]^{1/2}$}

\EndProcedure
\\ \\
{\it Note the use here of notation $\cG_{\btheta}(\bv)$ refers to tracing out the geodesic with initial velocity $\bv$ starting from point $\btheta$ (by mapping to the origin, using equation~\eqref{MVNgeodesic_solution_matrix}, and then mapping back). The base notation $\cG(\bv)$ is equivalent to $\cG_{(\b0,\bI)}(\bv)$.}

\end{algorithmic}
\caption{Adaptive path refinement for finding the geodesic between $\btheta$ and $\btheta'$}
\label{suppalg:MVNmultishooting_adaptive}
\end{algorithm}

\begin{algorithm}
\begin{algorithmic}
\Procedure{innerProduct}{$\bu, \bu', \btheta$}    \Comment{Inner product between two vectors at point $\btheta$}
\State \Return $\bumu^T \btheta_{\bSigma}^{-1} \bu'_{\bmu} + \frac{1}{2} \trop{ \btheta_{\bSigma}^{-1} \buSigma \btheta_{\bSigma}^{-1} \bu'_{\bSigma} }$    
\EndProcedure
\\
\Procedure{symKL}{$\btheta, \btheta'$}    \Comment{Evaluates symmeterised KL divergence between two MVNs}
\State \Return $\frac{1}{4} ( \btheta'_{\bmu} - \btheta_{\bmu} )^T \left( \btheta_{\bSigma}^{-1} + {\btheta'}_{\bSigma}^{-1} \right) ( \btheta'_{\bmu} - \btheta_{\bmu} ) + \frac{1}{4} \tr \left( \btheta_{\bSigma}^{-1} \btheta'_{\bSigma} + {\btheta'}_{\bSigma}^{-1} \btheta_{\bSigma} - 2\bI \right)$
\EndProcedure
\\
\Procedure{residualVector}{$\btheta, \btheta'$}    \Comment{Returns a vector joining $\btheta$ to $\btheta'$ in some sense}
\State $\btheta_t \gets \left(\btheta_{\bSigma}^{-1/2} (\btheta'_{\bmu} - \btheta_{\bmu}),\quad \btheta_{\bSigma}^{-1/2} \btheta'_{\bSigma} \btheta_{\bSigma}^{-1/2}\right)$    \Comment{Target point after transforming to origin}
\State \texttt{method} $\gets$ \texttt{\{euclid, taylor, eigen, projection\}}
\Switch{\texttt{method}}
    \Case{\texttt{euclid}}
        \State $(\bx, \bB) \gets \Bigl( (\btheta_t)_{\bmu}, \quad (\btheta_t)_{\bSigma} - \bI \Bigr)$     \Comment{Euclidean residual between origin and $\btheta_t$}
    \EndCase
    \Case{\texttt{taylor}}
        \State $(\bx, \bB) \gets $ Evaluate equation~\eqref{taylor_approximation}
    \EndCase
    \Case{\texttt{eigen}}
        \State $(\bx, \bB) \gets $ Evaluate equation~\eqref{MVNgeodesics_approx}
    \EndCase
    \Case{\texttt{projection}}
        \State $(\bx, \bB) \gets $ Match components in equation~\eqref{projection_approximation}
    \EndCase
    \\
\EndSwitch
\State $\br \gets \left(\btheta_{\bSigma}^{1/2} \bx, \quad \btheta_{\bSigma}^{1/2} \, \bB \, \btheta_{\bSigma}^{1/2}\right)$ \Comment{Transform the vector back into the original space}
\State \Return $\br$
\EndProcedure
\end{algorithmic}
\caption{Subfunctions used by the geodesic shooting routines}
\label{suppalg:subfunctions1}
\end{algorithm}

\begin{algorithm}
\begin{algorithmic}
\Procedure{backTransport}{$\br, \bv$}     \Comment{Transports $\br$ along geodesic defined by $\bv$}
\State \texttt{steps\_per\_length} $\gets 250$    \Comment{Specifies timestepping resolution}
\State $L \gets \Bigl[ \Call{innerProduct}{\bv, \bv, \bI} \Bigr]^{1/2}$     \Comment{Find length of geodesic}
\State $N \gets \Call{ceil}{\texttt{steps\_per\_length} * L}$    \Comment{Select number of numerical integration points}
\State $h \gets -1 / N$        \Comment{Timestep for geodesic numeric integration (backwards)}
\State $t = 1$       \Comment{Initialise numerical integration at geodesic endpoint}
\For{$i = 1, \ldots, N$}    \Comment{Loop from $t = 1$ back to $t = 0$}
\State $\br \gets \Call{RK4Step}{\textsc{parallelRHS}(\br, t, \bv), h, t}$    \Comment{RK4 for given RHS function}
\State $t \gets t + h$
\EndFor
\State \Return $\br$
\EndProcedure
\\
\Procedure{parallelRHS}{$\bu, t, \bv$}     \Comment{RHS of equations \eqref{parallel_transport} for updating $\bu$}
\State $\bg \gets \cG(t \bv)$     \Comment{Find geodesic endpoint at this time}
\State $\de{\bmu}{t} \gets \bg_{\bSigma} \bvmu$    \Comment{Find geodesic rate of change from equation~\eqref{MVNgeodesic_odes}}
\State $\de{\bSigma}{t} \gets \bg_{\bSigma} \left( \bvSigma - \bvmu \bp_{\bmu}^T  \right)$    \Comment{Find geodesic rate of change from equation~\eqref{MVNgeodesic_odes}}
\State \Return $\Biggl(  \frac{1}{2} \Bigl( \de{\bSigma}{t} \bg_{\bSigma}^{-1} \bu_{\bmu} + \bu_{\bSigma} \bg_{\bSigma}^{-1} \de{\bmu}{t} \Bigr), \quad    
\frac{1}{2} \Bigl( \de{\bSigma}{t} \bg_{\bSigma}^{-1} \bu_{\bSigma} + \bu_{\bSigma} \bg_{\bSigma}^{-1} \de{\bSigma}{t} - \de{\bmu}{t} \bu_{\bmu}^T - \bu_{\bmu} \left(\frac{d \bmu}{d t}\right)^T \Bigr) \Biggr) $
\EndProcedure
\end{algorithmic}
\caption{Additional subfunctions used by the geodesic shooting routines}
\label{suppalg:subfunctions2}
\end{algorithm}

\clearpage

\subsection{Supplementary Figures}

\begin{figure}[!h]
    \centering
    \includegraphics[width=0.75\textwidth, trim={1cm, 0cm, 1cm, 1.5cm}, clip]{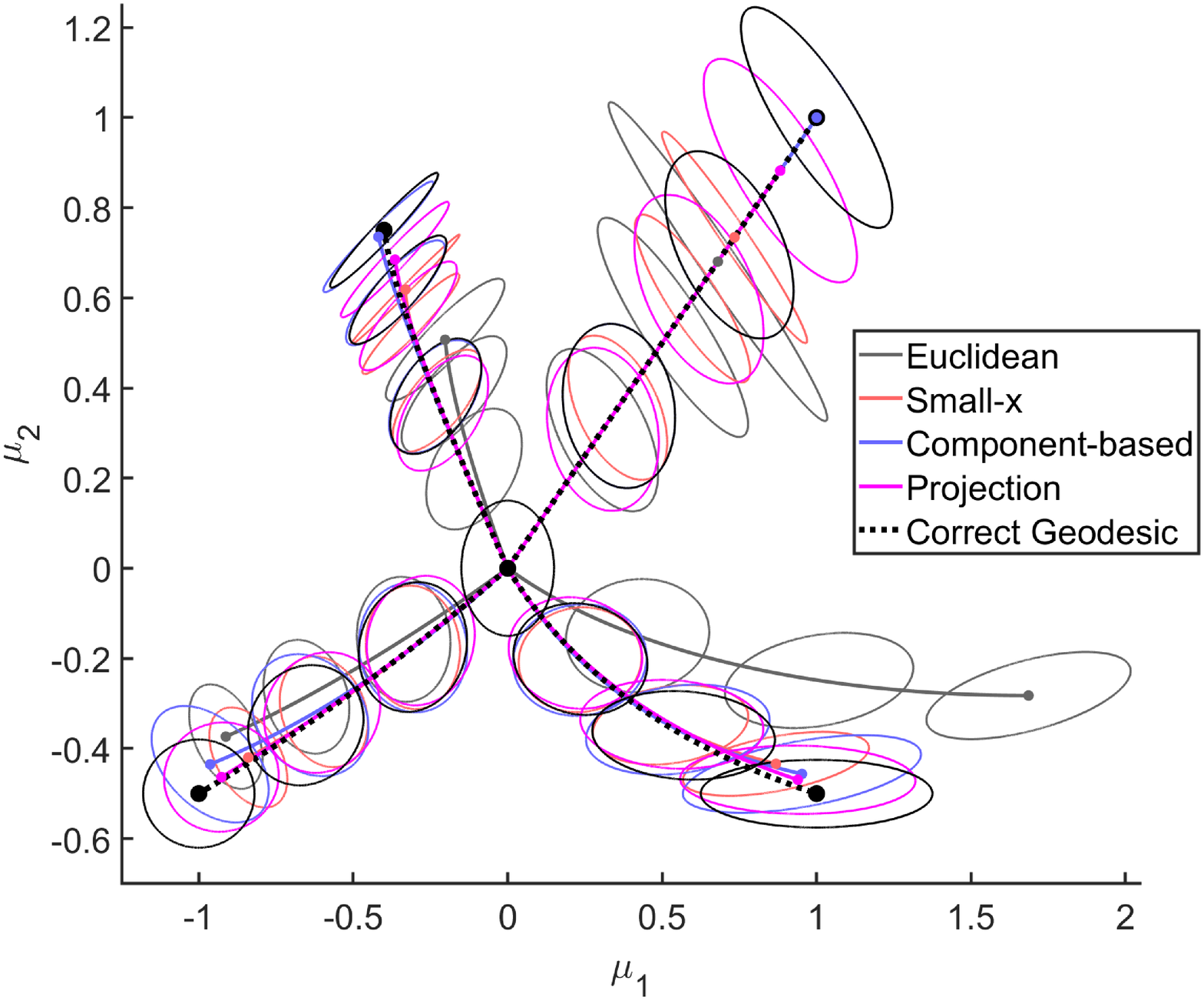}
    \caption{Visualisation of some example geodesics for bivariate MVNs, and the geodesics found by approximating the initial velocity to the target. The component-based approximation is exact for straight-line geodesics (top right), but the projection-based approximation is the best overall performer (see also \figref{suppfig:approximate_geodesic_performance}).}
    \label{suppfig:approximate_geodesic_examples}
\end{figure}

\begin{figure}
    \centering
    {\hspace{0.75cm} \large \bf{d = 2}  \hspace{6.25cm} \bf{d = 10}}\\
    \includegraphics[width=0.45\textwidth, trim={0.9cm, 0.1cm, 2.4cm, 1.25cm}, clip]{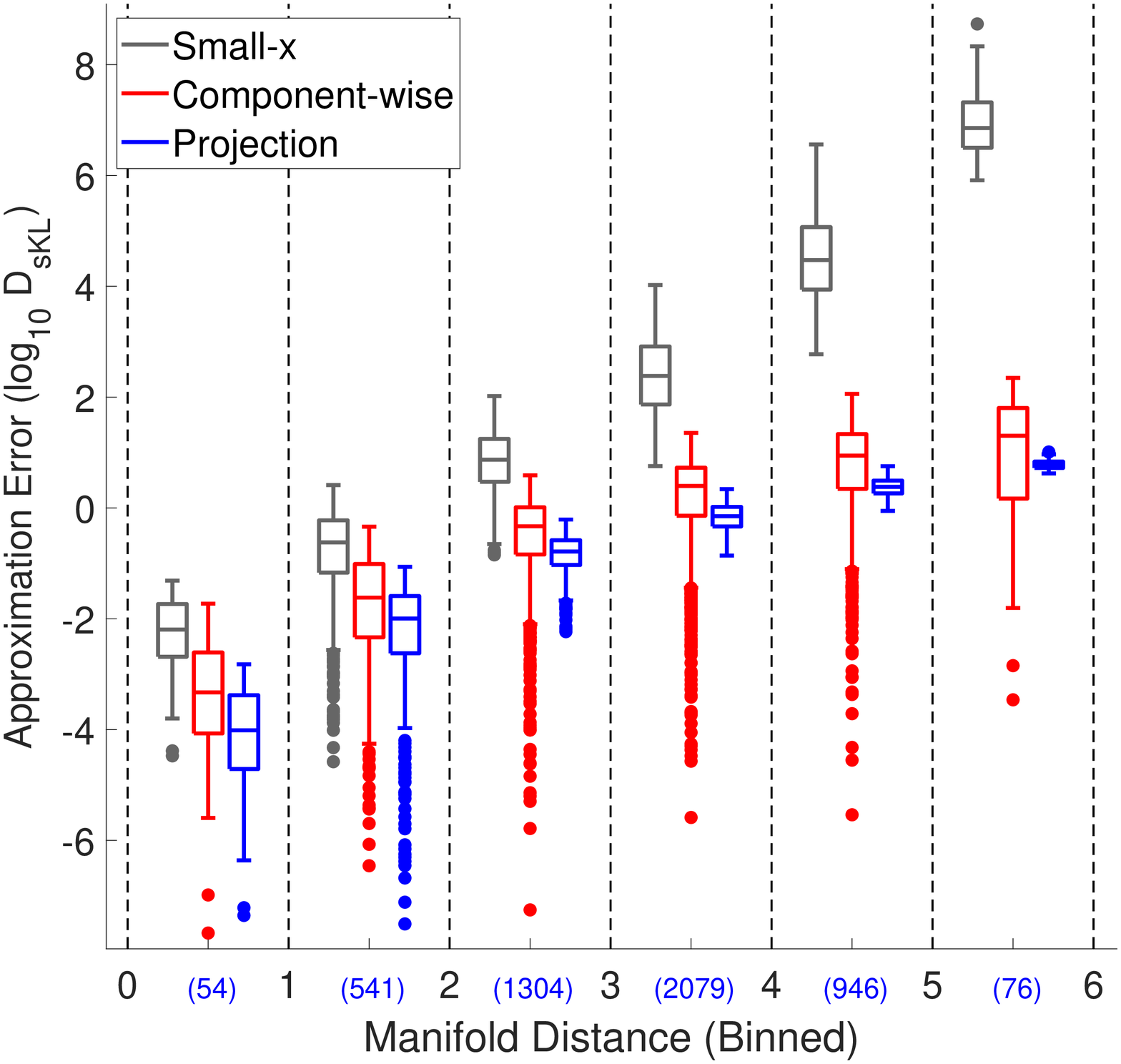}
    \hspace{0.05\textwidth}
    \includegraphics[width=0.45\textwidth, trim={0.9cm, 0.1cm, 2.4cm, 1.25cm}, clip]{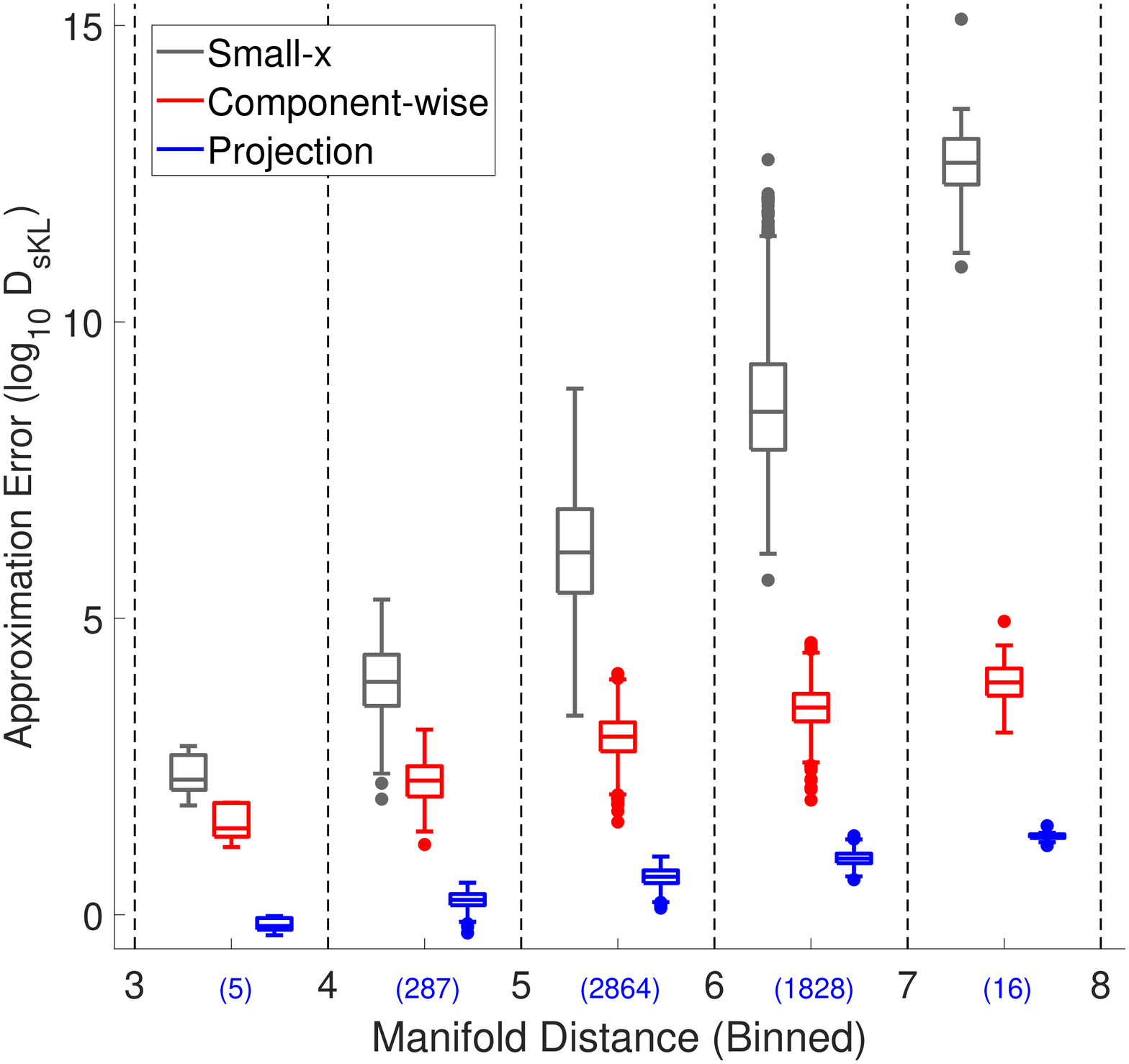}
    \caption{Performance (in terms of symmetric Kullback--Leibler divergence, see main text) of approximate geodesic velocities connecting the origin to 5000 randomly-generated target MVNs ($\bmu$ components selected $\mu_j \sim U[-5,5]$ and covariance eigenvalues $\log \lambda_j \sim U[-5, 5]$), for bivariate and 10-variate MVNs. Lengths are binned to aid in summarising the results (numbers in parantheses denote bin sample counts). Manifold projection performs best overall, and is noticeably more consistent for test problems generated in this fashion. However in the lower dimension, the component-wise approximation is sometimes superior for distant targets.}
    \label{suppfig:approximate_geodesic_performance}
\end{figure}

\begin{figure}
    \centering
    \includegraphics[width=0.65\textwidth,trim={0cm, 0cm, 2cm, 0cm},clip]{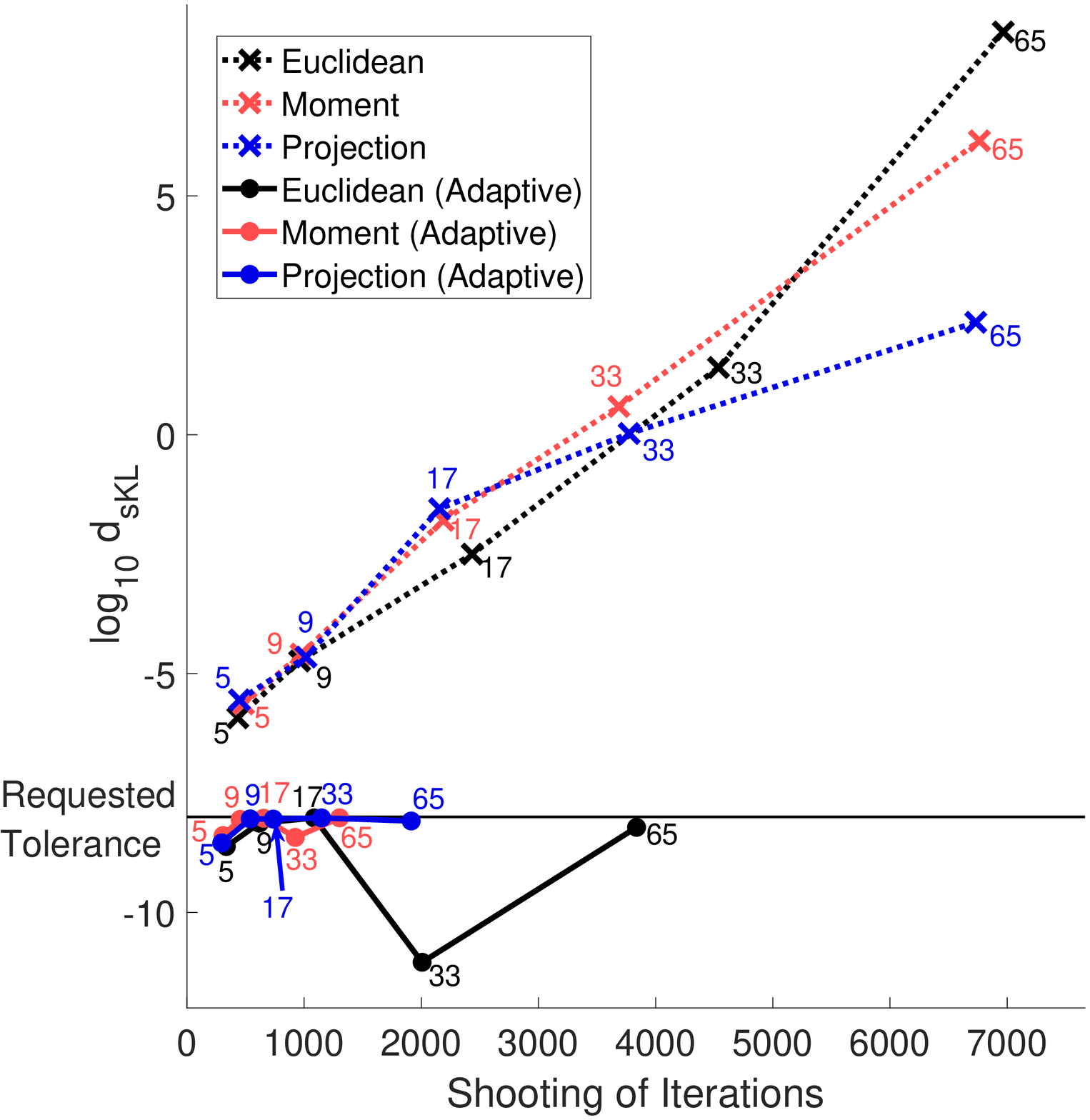}
    \caption{Standard path refinement (dashed lines) fails to find a geodesic to the tolerance used for shooting subproblems, and suffers greatly from choosing too many path points. Adaptive path refinement (solid lines) rectifies these issues, but still benefits from choosing as few points as possible. Lines are drawn to aid figure interpretation only, and simply connect the points of data for each method.}
    \label{suppfig:path_refine_Npoints}
\end{figure}

\end{document}